# Tutorial and Practice in Linear Programming

## Optimization Problems in Supply Chain and Transport Logistics

Raj Bridgelall, Ph.D.
Associate Professor
Department of Transportation, Logistics, and Finance, North Dakota State University
Phone: (408) 607-3214; Email: raj@bridgelall.com
ORCID: 0000-0003-3743-6652

**Abstract**

This tutorial is an andragogical guide for students and practitioners seeking to understand the fundamentals and practice of linear programming. The exercises demonstrate how to solve classical optimization problems with an emphasis on spatial analysis in supply chain management and transport logistics. All exercises display the Python programs and optimization libraries used to solve them. The first chapter introduces key concepts in linear programming and contributes a new cognitive framework to help students and practitioners set up each optimization problem. The cognitive framework organizes the decision variables, constraints, the objective function, and variable bounds in a format for direct application to optimization software. The second chapter introduces two types of mobility optimization problems (shortest path in a network and minimum cost tour) in the context of delivery and service planning logistics. The third chapter introduces four types of spatial optimization problems (neighborhood coverage, flow capturing, zone heterogeneity, service coverage) and contributes a workflow to visualize the optimized solutions in maps. The workflow creates decision variables from maps by using the free geographic information systems (GIS) programs QGIS and GeoDA. The fourth chapter introduces three types of spatial logistical problems (spatial distribution, flow maximization, warehouse location optimization) and demonstrates how to scale the cognitive framework in software to reach solutions. The final chapter summarizes lessons learned and provides insights about how students and practitioners can modify the Phyton programs and GIS workflows to solve their own optimization problem and visualize the results.

# LIST OF FIGURES

# LIST OF TABLES

# Chapter 1 - Introduction

## 1.1 Background

Some of the optimization topics covered in this tutorial will be familiar to those who studied operations research (OR), management science, or decision science. Students can use this tutorial to learn the fundamentals and practice of solving optimization problems in the broader fields of supply chain management and transport logistics.

## 1.2 Concepts in Linear Programming

The term *linear programming* arises from the fact that the objective function is a *linear* combination of *decision variables* and *parameters* that one seeks to maximize or minimize. For example, classic problems seek to maximize profits and flow and to minimize cost or time. The parameters in the linear combination of variables are fixed values that represent facts of the problem. A program solves the problem based on constraints of the variables or outputs. The model *constraints* are also linear combinations of the decision variables and the resource or input parameters. The value of constraints are also parameters that represent facts of the problem. The *intersection* of all constraints forms a *feasible solution space* such that a unique combination of variables will result in the *optimum solution*. The name, *decision variables,* suggests that the optimum values to achieve the objective could guide decision making. Some optimization problems are *infeasible*. Even if an optimization problem is *feasible*, it may have no optimum solution and, therefore, becomes an *unbounded* problem.

## 1.3 Types of Linear Programming

Linear programming can be integer linear programming (ILP), binary integer programming (BIP), and mixed integer linear programming (MILP). The decision variables of ILP are positive integers, including zero. The decision variables of BIP are binary—they represent Boolean logic by assigning values of zero or one. The decision variables of MILP can be a mix of continuous, integer, and binary variables.

## 1.4 Problem Formulation Equations

The following are steps to formulate the optimization problem:
1) Define a set of **decision variables**—these are the choices that someone must make to achieve the optimal solution under various constraints on those choices.
2) Define the **constraints** as linear functions of the decision variables with parameters defined to represent a model of the system.
3) Define the **objective function** as a *single* linear function of the decision variables with parameters that represent one unit of the associated decision variable. Examples of parameters are cost in dollar units and flow in units of vehicles per hour.
4) Define the lower and upper **bounds** of each decision variable. Binary decision variables can take on finite values of zero or one and positive variables have a lower bound of zero.

The **objective** value $Z$ is the dot product of the objective coefficients and the decision variables where

$$Z = \sum_{i=1}^{N} c_i X_i \qquad (1)$$

is the objective function or the *functional relationship*. The optimization problem is to determine the values of $X_i$ that will either maximize or minimize Z, subject to the constraints. The **constraints** can be bounds. For example, in a production optimization problem, a bound on each decision variable can be the number of items $i$ produced where

$$L_i \leq X_i \leq U_i \qquad (2)$$

$L_i$ and $U_i$ are the lower and upper bounds, respectively, for item $i$. Additional constraints can establish the total of any items or combination of items produced. This tutorial introduces a cognitive framework to help students and practitioners formulate a problem by using a standard structure dubbed the problem formulation table (PFT).

## 1.5 Problem Formulation Table (PFT)

Table 1 shows a general structure for the proposed cognitive framework to help organize facts of the problem such as the decision variables, constraint parameters, objective function parameters, and bounds of the decision variables. This PFT is scalable so that a program can read the table as several matrices and convert the problem statement into a *general form* required by the selected optimizer program.

Table 1: Variables and Constants associated with an Optimization Problem

| Variables (Outputs) | Constraints | | | | Objective (Min/Max) | Variable Bounds | |
|---|---|---|---|---|---|---|---|
| | $J_1$ | $J_2$ | ... | $J_K$ | | L | U |
| $X_1$ | $\alpha_{11}$ | $\alpha_{12}$ | | $\alpha_{1K}$ | $c_1$ | $L_1$ | $U_1$ |
| $X_2$ | $\alpha_{22}$ | $\alpha_{22}$ | | $\alpha_{2K}$ | $c_2$ | $L_2$ | $U_2$ |
| ⋮ | ⋮ | ⋮ | ⋮ | ⋮ | ⋮ | ⋮ | ⋮ |
| $X_N$ | $\alpha_{N1}$ | $\alpha_{N2}$ | | $\alpha_{NK}$ | $c_N$ | $L_N$ | $U_N$ |
| | $b_1$ | $b_2$ | | $b_K$ | Z | | |

The decision variables $X_i$ are the individual outputs that the optimization program needs to find. In a production problem, the $N$ decision variables might be individual products that a factory must produce. The decision variables have indices $i = 1, 2, ..., N$. The $K$ constraints can be the various resources, parts, or ingredients needed to build each product and they have indices $j = 1, 2, ..., K$. The fixed *parameters* values $\alpha_{ij}$ represent the proportion of resource $J_j$ needed to produce a single item $X_i$. For example, $\alpha_{ij}$ could be the number of hours or fraction of hours needed to produce a widget or deliver a service $X_i$. The fixed values $c_i$ are proportions of the objective unit associated with each decision variable $X_i$. For example, $c_i$ could be the cost to produce a single widget $X_i$ or the profit earned from selling a single widget $X_i$. The *objective coefficients* become the coefficients of the *objective function*. The last two columns of the PFT contain the lower (L) and upper (U) bounds of each decision variable. The PTF need not specify these columns when all decision variables are binary because they have the same bounds.

 Constraints on resource $j$ is the dot product of the resource parameter vector $J_j$ and the decision variable vector $X$ with an *upper bound* given by the value $b_j$ in the last row. That is, the standard formulation, as required by the optimization program, would be

$$\sum_{i=1}^{N} \alpha_{ij} X_i \leq b_j \tag{3}$$

where $b_j$ is the maximum amount or available capacity of resource $j$. Constraints for one decision variable can be a function of other decision variables. For example, if for a constraint $j$ the number of item 2 should be twice the number of item one, then constraint $j$ will be

$$X_2 = 2X_1 \tag{4}$$

which when converted to the *standard form* becomes

$$-2X_1 + X_2 = 0. \tag{5}$$

Hence the parameter values would be $\alpha_{1j} = -2$, $\alpha_{2j} = 1$, $\alpha_{ij} = 0$, and $b_j = 0$.

## 1.6 PFT Evaluation

The model builder should evaluate the PFT to identify trivial infeasibilities, trivial unboundedness, and any issues in the problem formulation by checking for the following:
1. Any column containing all zeros represents a trivial constraint.
2. A column singleton (a single 1) in an equality constraint matrix represents a fixed value for a decision variable. Hence, the model should not include that value in the optimization because it is a *constant*, not a *variable*.
3. Any row containing all zeros represents an unconstrained variable.
4. A row singleton (a single 1) in an inequality constraint matrix represents a simple bound and the model builder could use it to simplify the problem formulation.

## 1.7 Software Tools

### 1.7.1 Python Program

Python is popular in data science because of its access to a diversity of standard and specialized libraries. Python is an interpreted rather than a compiled language. This means that coding results are instantaneous because there is no need to compile code before executing it. Python is suitable for any scale project, and it is portable across many operating systems. Downloading and running the Anaconda distribution is one of the easiest ways to install Phyton (Anaconda Inc., 2022). This tutorial recommends using the Spyder IDE to execute the Python code examples presented. Copy and paste the source code from this document directly into the Spyder IDE code window. For all the programs in this guide, change the "datapath_in =" variable to match the location of those files on your computer. Also, ensure that the "infile" variable matches the names of the data files (Excel or CSV) that you produce when copying and pasting data tables.

### 1.7.2 MIP Library

Numerous Python libraries facilitate mixed integer programming. This tutorial selected the MIP library from COIN-OR (Santos & Toffolo, 2020). The tools feature high level modeling with rapid code execution. It is compatible with Pypy, which is a compiler that can make large MIP programs run many times faster. Install the MIP library by launching the Anaconda prompt and using the PIP command as suggested in the COIN-OR documentation.

### *1.7.3 QGIS®*

QGIS is a free open-source GIS tool that runs on many different types of operating systems, including Windows®, iOS®, and Linux (QGIS, 2022). The standard installation offers many of the same features of popular commercial software, and many more through free open-source plugins. Follow the instructions on the QGIS website to download and install the latest stable version.

### *1.7.4 GeoDA™*

GeoDA™ is free software developed by a team at the University of Chicago, Center for Spatial Data Science, and partially funded by a grant from the National Science Foundation (GeoDA, 2022). GeoDA focuses on tools to conduct a large variety of spatial analysis and has basic GIS methods necessary to perform statistical operations that involve points and polygon data. As of 2022, the tool does not support line data. Therefore, it is a complement to QGIS. One of the most important features of GeoDA is linking, which is an ability to highlight selected data across all analysis windows. Another feature called *brushing* is the ability to move a selection tool in any window to dynamically highlight the selected data across maps and charts. Linking and brushing are powerful tools in exploratory spatial data analysis (ESDA) that can enable analyst to observe relationships across various visualizations of the data to gain insights not otherwise possible. For example, using a box plot tool to select the outliers will highlight their locations on a map, as well as their attributes in the data table. Download and install the latest *stable* version of GeoDA from their website.

## 1.8  Recommended Reading

1. Shaw, Shih-Lung. "Geographic information systems for transportation–An introduction." *Journal of Transport Geography* 3, no. 19 (2011): 377-378. (Shaw, 2011)

2. Miller, Harvey J., and Shih-Lung Shaw. "Geographic information systems for transportation in the 21st century." *Geography Compass* 9, no. 4 (2015): 180-189. (Miller & Shaw, 2015)

3. Anselin, Luc, Ibnu Syabri, and Youngihn Kho. "GeoDa: an introduction to spatial data analysis." *Geographical analysis* 38, no. 1 (2006): 5-22 [pdf]. (Anselin, Syabri, & Kho, GeoDa: an introduction to spatial data analysis, 2010)

4. Anselin, Luc. "The GeoDA Book: Exploring Spatial Data." Chapter 1, GeoDa Press LLC, Chicago, IL, 2017 [pdf]. (Anselin, 2017 )

# Chapter 2 - Mobility Optimization

This chapter describes how to apply the PFT introduced in Chapter 1 to set up and solve two popular optimization problems in mobility optimization: Shortest Path and Minimum Cost Tour.

## 2.1 Shortest Path

The goal of the shortest path problem is to identify a set of links that form a path connecting a starting node *s* with a terminal node *t* such that the cost in terms of travel distance or travel time is minimum. Table 2 shows how to use a PFT to organize the optimization problem.

Table 2: PFT for the Shortest Path Problem

| Links | Node *j* Flow Constraint | | | | Node *j* Input Constraint | | | | Travel Time (Hours) |
|---|---|---|---|---|---|---|---|---|---|
| | 1 | 2 | … | K | 1 | 2 | … | K | |
| $X_{11}$ | $\alpha_{11}$ | $\alpha_{12}$ | | $\alpha_{1K}$ | $\beta_{11}$ | $\beta_{12}$ | | $\beta_{1K}$ | $c_1$ |
| $X_{21}$ | $\alpha_{21}$ | $\alpha_{22}$ | | $\alpha_{2K}$ | $\beta_{21}$ | $\beta_{22}$ | | $\beta_{2K}$ | $c_2$ |
| ⋮ | ⋮ | ⋮ | ⋮ | ⋮ | ⋮ | ⋮ | ⋮ | ⋮ | ⋮ |
| $X_{KK}$ | $\alpha_{N1}$ | $\alpha_{N2}$ | | $\alpha_{NK}$ | $\beta_{N1}$ | $\beta_{N2}$ | | $\beta_{NK}$ | $c_N$ |
| | -1 | 0 | | 1 | 1 | 1 | | 1 | T |

Each decision variable $X_{ij}$ represents a single link that exists between node *i* and node *j*. The optimization selects a link $X_{ij}$ from a set {***L***} of *N* defined links by assigning the variable a value of 1, and otherwise assigns a value of 0. The label of the starting node *s* is 1 and the label of the terminal node *t* is the total number of nodes *K* in the network. To represent the network digitally, the analyst establishes a set of flow constraints by assigning $\alpha_{ij}$ one value from the set {0, 1, -1}. The flow constraint of each node *j* assigns $\alpha_{ij} = 1$ for all links entering the node, -1 for all links exiting the node, and 0 for all other values of $\alpha_{ij}$. The *exit* constraint for the starting node *s* is $b_1 = -1$ because exactly one unit of flow must exit that node. Similarly, the *entry* constraint for the terminal node is $b_K = +1$ because the same single unit of the flow must enter the *t* node. A *flow conservation* constraint requires the assignment of $b_j = 0$ ($j \neq s, j \neq t$) to all other intermediate nodes because their total units of flow entering must also exit them. Another *entry* constraint assures that the unit flow must enter a node from one and only one of the links connected to it. The notation $X_{ij}^k$ indicate the set of entry links connected to node *k*. The optimization problem becomes:

Minimize 
$$T = \sum_{X_{ij} \in \{L\},\ n=1}^{N} c_n X_{ij} \qquad (6)$$

Subject to:

$$\sum_{X_{ij}\in\{L\},\ n=1}^{N} \alpha_{ij}X_{ij} = \begin{cases} -1 & if\ i = s \\ +1 & if\ i = t \\ 0 & otherwise \end{cases} \quad (7)$$

and

$$\sum_{X_{ij}\in\{L\},\ k=1}^{K} X_{ij}^k \leq 1 \quad (8)$$

where

$$X_{ij} \in \{0, 1\} \quad (9)$$

For nodes not selected on the path, the sum of input flows will be zero, and for nodes that are selected on the path, the sum of input flows will be exactly one, hence the reason for the ≤ in the second set of constraints.

### 2.1.1 Example Problem

This example is adopted from a problem posed by Taylor (Bernard W. Taylor, 2019). The problem is to find the shortest trucking route between Los Angeles and St. Louis. Figure 1 shows a network representation of typical truck routes between Los Angeles and St. Louis.

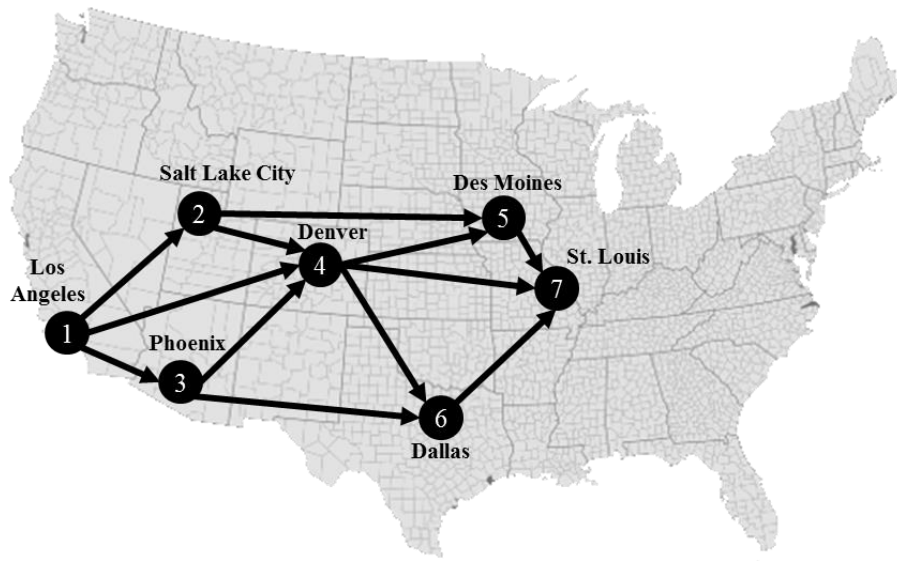

Figure 1: Network representation for the shortest shipping route problem.

Table 3 shows the PFT to solve the shortest path problem with an optimization software.

Table 3: PFT for the Shortest Path Problem Exercise

| Net Links | Node *j* Flow Constraints | | | | | | | Node *j* Input Constraints | | | | | | | Road Miles | Geodesic Miles |
|---|---|---|---|---|---|---|---|---|---|---|---|---|---|---|---|---|
| | *1* | *2* | *3* | *4* | *5* | *6* | *7* | *1* | *2* | *3* | *4* | *5* | *6* | *7* | | |
| $X_{12}$ | -1 | 1 | | | | | | | 1 | | | | | | 688 | 579.4 |
| $X_{13}$ | -1 | | 1 | | | | | | | 1 | | | | | 373 | 357.6 |
| $X_{14}$ | -1 | | | 1 | | | | | | | 1 | | | | 1016 | 831.3 |
| $X_{24}$ | | -1 | | 1 | | | | | | | 1 | | | | 519 | 371.8 |
| $X_{25}$ | | -1 | | | 1 | | | | | | | 1 | | | 1066 | 953.1 |
| $X_{34}$ | | | -1 | 1 | | | | | | | 1 | | | | 821 | 585.7 |
| $X_{36}$ | | | -1 | | | 1 | | | | | | | 1 | | 1065 | 885.8 |
| $X_{45}$ | | | | -1 | 1 | | | | | | | 1 | | | 671 | 610.8 |
| $X_{46}$ | | | | -1 | | 1 | | | | | | | 1 | | 794 | 661.7 |
| $X_{47}$ | | | | -1 | | | 1 | | | | | | | 1 | 851 | 796.7 |
| $X_{57}$ | | | | | -1 | | 1 | | | | | | | 1 | 349 | 273.3 |
| $X_{67}$ | | | | | | -1 | 1 | | | | | | | 1 | 631 | 547.8 |
| =, = | -1 | 0 | 0 | 0 | 0 | 0 | 1 | 1 | 1 | 1 | 1 | 1 | 1 | 1 | R | G |

The nodes represent the location of major cities, and the arcs represent established travel routes between the nodes. The first column of the PFT contains the decision variables, which are also the labels for each arc. The last two columns contain a value associated with each decision variable, which in this case is a distance in miles to compute the minimum total distance. The road miles indicated is based on the recommended route using Google Maps®. The Geodesic miles indicated is the geodesic distance derived using a GIS software. For illustration, the trip ends are truck stop locations at the approximate center of each city. The objective function is the *dot product* of the decision variables and one of the cost columns. Starting with the road route costs, the optimization problem is:

Minimize:

$$R = 688X_{12} + 373X_{13} + 1016X_{14} + 519X_{24} + 1066X_{25} + 821X_{34} \\ + 1065X_{36} + 671X_{45} + 794X_{46} + 851X_{47} + 349X_{57} + 631X_{67} \quad (10)$$

The seven node flow constraint columns define the network structure by defining the entry and exit arc flow for each node. A -1 is an exit flow and a +1 is an entry flow. The last row indicates the constraint constant for the *dot product* of its respective column. In this example, the dot product must be equal to the constant. This results in the following constraints:

Subject to:
$$-X_{12} - X_{13} - X_{14} = -1$$
$$\text{or} \tag{11}$$
$$X_{12} + X_{13} + X_{14} = 1$$
$$X_{12} - X_{24} - X_{25} = 0 \tag{12}$$
$$X_{13} - X_{34} - X_{36} = 0 \tag{13}$$
$$X_{14} + X_{24} + X_{34} - X_{45} - X_{46} - X_{47} = 0 \tag{14}$$
$$X_{36} + X_{46} - X_{57} = 0 \tag{15}$$
$$X_{36} + X_{46} - X_{67} = 0 \tag{16}$$
$$X_{47} + X_{57} + X_{67} = 1 \tag{17}$$
where:
$$X_{ij} \in \{0, 1\} \tag{18}$$

### 2.1.2 Solution Exercise

The optimizer from the MIP library solves linear programming problems in the following form:
$$\min_{x} \boldsymbol{c}^T \boldsymbol{x} \tag{19}$$
such that
$$\boldsymbol{A}^T \boldsymbol{x} = \boldsymbol{b} \tag{20}$$
where
$$\boldsymbol{x} \in \{0, 1\} \tag{21}$$

The vector $\boldsymbol{x}$ are the binary decision variables, $\boldsymbol{c}$ is a vector of objective function coefficients, $\boldsymbol{A}$ is the constraint matrix extracted from the PFT, and $\boldsymbol{b}$ is a vector of the constants on the right side of the *equality* or *inequality* constraint equations. That is, each row of $\boldsymbol{A}^T$ contains the coefficients of a linear equality or inequality constraint on the decision variables in vector $\boldsymbol{x}$. The solution corresponds to the defined order of the decision variables.

### 2.1.3 Exercise on BIP with Python

Create a copy of the PFT (Table 3) in an Excel file. Modify and run the Python program below. Change the variable "datapath_in" to match the path where the Excel file is located on your computer. Change the "infile" variable to match the name of the Excel file on your computer. The program extracts the PFT by skipping the first row to create a data frame. Because Python imports blank cells as "nan" the program replaces them with zeros. The program stores the variable names in the first column as a series of strings. Per equation (11), the program then converts the first constraint column to positive values so that the values are binary. The program extracts the parameters into matrices and arrays that meets the MIP requirements for adding variables, objectives, and constraints to the model object.

```python
# Author: Dr. Raj Bridgelall (raj.bridgelall@ndsu.edu)
# Shortest Path MIP
from IPython import get_ipython
get_ipython().magic('clear')                            # Clear the console
import pandas as pd
from pathlib import Path
from mip import *                        # install library from Anaconda prompt: pip install mip
import numpy as np
#%%
datapath_in = 'C:/Users/Admin/Documents/OneDrive/Documents/UGPTI/Teaching/TL 885/Lectures/6 - Mobility Optimization/Lab/'
infile = 'Shortest Path PFT.xlsx'                       # Input filename
filepath_in = Path(datapath_in + infile)                # Path name for untruncated signal
df = pd.DataFrame(pd.read_excel(filepath_in, skiprows=1))    # Read Problem Formulation Table
#%%
Nc = int(np.floor(df.shape[1]/2)-1)     # Number of resource constraints
Nd = df.shape[0]-1                      # Number of decision variables
#%%
df = df.fillna(0)                 # replace all NaN (blanks in Excel) with zeros
VarName = df.iloc[0:Nd,0]         # String of variable names
df[1] = df[1].astype(int) * -1    # Convert starting node column to positive constraint of +1
A1_parms = df.iloc[0:Nd,1:Nc+1].astype(int)             # Extract the constraints cols
b1_parms = list(df.values[-1][1:Nc+1].astype(int))      # Get constraint parameters in an array
A2_parms = df.iloc[0:Nd,Nc+2:-2].astype(int)            # Extract the constraints cols (except 1st)
b2_parms = list(df.values[-1][Nc+2:-2].astype(int))     # Get constraint parameters in an array
c_parms = df.values.transpose()[-2][:-1].astype(int)    # Get objective parameters in an array
#%%
m = Model(solver_name=CBC)                              # Instantiate the solver model object
x = [m.add_var(name=VarName[i],var_type=BINARY) for i in range(Nd)] # define decision vars in x
m.objective = minimize(xsum(c_parms[i]*x[i] for i in range(Nd)))    # add objective function
#%% Add each constraint column
for j in range(Nc):
    m.add_constr( xsum(A1_parms.iloc[i,j]*x[i] for i in range(Nd)) == b1_parms[j],\
              "Cons1_"+str(j+1) )

for j in range(Nc-1):
    m.add_constr( xsum(A2_parms.iloc[i,j]*x[i] for i in range(Nd)) <= b2_parms[j],\
              "Cons2_"+str(j+2) )
#%%
Status = m.optimize()
#%% Print the results
print('Model has {} vars, {} constraints and {} nzs'.format(m.num_cols, m.num_rows, m.num_nz))
selected = [VarName[i] for i in range(Nd) if x[i].x != 0]
print('Selected Paths: {} '.format(selected))
print('Minimum Cost = {} Miles'.format(m.objective_value))
for j in range(Nc+Nc-1):
    print(m.constrs[j])
print("Number of Solutions = ", m.num_solutions)
print("Status = ", Status)
```

The program prints the solution as follows:
```
Model has 12 vars, 13 constraints and 36 nzs
Selected Paths: ['X14', 'X47']
Minimum Cost = 1867.0 Miles
Cons1_0: +1.0 X12 +1.0 X13 +1.0 X14 = 1.0
Cons1_1: +1.0 X12 -1.0 X24 -1.0 X25 = -0.0
Cons1_2: +1.0 X13 -1.0 X34 -1.0 X36 = -0.0
Cons1_3: +1.0 X14 +1.0 X24 +1.0 X34 -1.0 X45 -1.0 X46 -1.0 X47 = -0.0
Cons1_4: +1.0 X25 +1.0 X45 -1.0 X57 = -0.0
Cons1_5: +1.0 X36 +1.0 X46 -1.0 X67 = -0.0
Cons1_6: +1.0 X47 +1.0 X57 +1.0 X67 = 1.0
Cons2_0: +1.0 X12 <= 1.0
Cons2_1: +1.0 X13 <= 1.0
Cons2_2: +1.0 X14 +1.0 X24 +1.0 X34 <= 1.0
Cons2_3: +1.0 X25 +1.0 X45 <= 1.0
Cons2_4: +1.0 X36 +1.0 X46 <= 1.0
Cons2_5: +1.0 X47 +1.0 X57 +1.0 X67 <= 1.0
Number of Solutions =  1
Status =  OptimizationStatus.OPTIMAL
```

Note that the input constraints for node $j = 1$ is all zeros, so the program eliminates that constraint when building the model. The program then calls the optimize function, which returns a status that indicates the results. The output indicates that the optimum solution is to take routes $X_{14}$ and $X_{47}$, which agrees with a visual observation that it should be the shortest path. This solution results in a total of 1,867 miles along primary roads. The equivalent cost in fuel consumed is the total miles divided by the average miles per gallon achieved by the truck.

### 2.1.4 Exercise with GIS

Use QGIS to determine the geodesic distances for the arcs between the cities. Those distances may represent more direct travel by air instead of by road. The course associated with this tutorial provides step-by-step instructions on how to complete this exercise. The last column of Table 3 shows the resulting values in miles. As an exercise, change the code to use those values as the new cost function, rerun the program, and compare the solution.

### 2.1.5 Further Reading

1. Taccari, Leonardo. "Integer programming formulations for the elementary shortest path problem." *European Journal of Operational Research* 252, no. 1 (2016): 122-130. (Taccari, 2016)

## 2.2 Minimum Cost Tour

The famous "traveling salesman" problem is a popular example in this category of problems. The problem is to determine how to visit all nodes on a network at most once while minimizing the travel cost of the tour. A node on the network represents a city that the salesman must visit. The value associated with each arc is the cost to travel between those cities. A tactic to generalize the problem of defining the network is to set an arc value to infinity where a direct connection between a pair of cities (points) does not exist. The mathematical formulation for the optimization problem is:

Minimize

$$D = \sum_{i=1}^{N} \sum_{j=1, j \neq i}^{N} c_{ij} X_{ij} \qquad (22)$$

Subject to:

$$\sum_{i=1, i \neq j}^{N} X_{ij} = 1 \qquad (j = 1, 2, \ldots, N) \qquad (23)$$

and

$$\sum_{j=1, j \neq i}^{N} X_{ij} = 1 \qquad (i = 1, 2, \ldots, N) \qquad (24)$$

and

$$u_i - u_j + 1 \leq (N - 1)(1 - X_{ij}) \qquad (\forall i \neq 1, \forall j \neq 1) \qquad (25)$$

where

$$X_{ij} = \begin{cases} 1 & \text{the path goes from node } i \text{ to node } j \\ 0 & \text{otherwise} \end{cases} \quad (26)$$

$$u_1 = 1, \quad 2 \leq u_i \leq N \quad (\forall i \neq 1) \quad (27)$$

The optimization is to minimize the total distance $D$ traveled on the tour when $c_{ij}$ represents the distance between cities $i$ and $j$. The first constraint assures that each city $j$ is arrived at from exactly one other city $j$. The second constraint assures that a departure from each city $i$ is to exactly one other city $j$. Hence, the first two constraints enforce a single-entry-single-exit (SESX) condition for all cities.

The SESX constraints are not enough to solve the full tour problem because the mathematical solution could result in one or more isolated sub-tours that meet those constraints for all cities. For example, Figure 2 shows a solution with subtours that meet the SESX constraint.

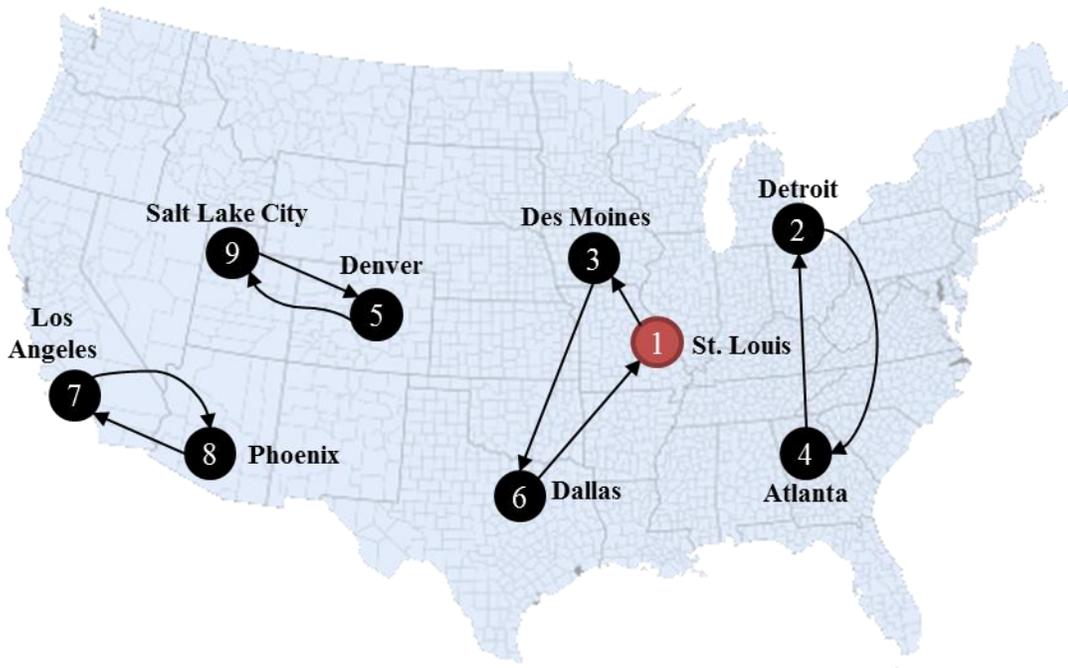

Figure 2: Subtour examples that satisfy the SESX condition for all cities.

The salesperson visited all cities exactly once but leaving any subtour to complete the full tour would violate the SESX constraint. Therefore, another constraint is needed to introduce an ordering of the visits without adding new decision variables. Adding "dummy variables" $u_i$ that are associated with each of the $N$ network nodes achieves an ordering constraint. The first dummy variable, associated with the starting node, gets a value of 1. The other dummy variables can take on integer values between 0 and $N$-1. The right side of the inequality becomes zero if the algorithm assigns arc $X_{ij} = 1$ because $(1-X_{ij}) = 0$. This assignment then forces $(u_i - u_j) = -1$ because the left side of the inequality must become zero (-1 + 1 = 0), and that forces $u_i$ to be exactly one position in the tour lower than $u_j$.

Table 4 shows the PFT for a minimum cost tour problem. There are N × (N-1) decision variables that represent the possible number of arcs between the cities. Note that when there are more than a few nodes involved, it is not practical to use a PFT. Here, the PFT serves more as a cognitive tool to organize the variables, constraints, and optimization. The PFT does not include decision variables for identity nodes where $i = j$ because those are not involved in a tour. Because $u_1 = 1$ is a trivial constraint, the table includes it as a node.

Table 4: PFT for the Minimum Cost Tour Problem Exercise

| All Links | Entry Constraints | | | | Exit Constraints | | | | Subtour Constraints (i, j) | | | | | Distance (Miles) |
|---|---|---|---|---|---|---|---|---|---|---|---|---|---|---|
| | $j_1$ | $j_2$ | ... | $j_N$ | $i_1$ | $i_2$ | ... | $i_N$ | 1,2 | ... | 1, N | N,1 | ... | N-1, N | |
| $X_{12}$ | 1 | | | | 1 | | | | N | | | | | | $c_{11}$ |
| ⋮ | | 1 | | | 1 | | | | | | | | | | ⋮ |
| $X_{1N}$ | | | 1 | | 1 | | | | | | | | | | $c_{1N}$ |
| $X_{21}$ | 1 | | | | | 1 | | | | | | | | | $c_{21}$ |
| ⋮ | | 1 | | | | 1 | | | | | | | | | ⋮ |
| $X_{2N}$ | | | 1 | | | 1 | | | | | | | | | $c_{2N}$ |
| ⋮ | | | | | | | | | | | | | | | ⋮ |
| $X_{N1}$ | 1 | | | | | | | 1 | | | | | | | $c_{N1}$ |
| ⋮ | | 1 | | | | | | 1 | | | | | | | ⋮ |
| $X_{N-1,N}$ | | | 1 | | | | | 1 | | | | | | N | $c_{N\_N-1}$ |
| $u_2$ | | | | | | | | | 1 | | | | | | |
| ⋮ | | | | | | | | | | | | | | 1 | |
| $u_N$ | | | | | | | | | | | | | | | |
| $u_2$ | | | | | | | | | -1 | | | | | | |
| ⋮ | | | | | | | | | | | | | | | |
| $u_N$ | | | | | | | | | | | | | | -1 | |
| =, =, ≤ | 1 | 1 | 1 | 1 | 1 | 1 | 1 | 1 | N-1 | | | | | N-1 | D |

Rewriting the subtour constraint in the standard MIP format yields

$$u_i - u_j + X_{ij}N \leq (N - 1) \qquad (\forall i \neq 1, \forall j \neq 1) \qquad (28)$$

such that the right side of the inequality is a constant equal to the number of nodes visited, excluding the home node. Note that the logic holds in that if $X_{ij} = 1$, then $u_i - u_j$ must be -1 so that the inequality holds. That is, if $X_{ij} = 1$, then the equation must evaluate to $-1 + N \leq (N - 1)$. Note that $u_i - u_j = -1$ is equivalent to $u_j = u_i + 1$. This forces the next node visited to have a dummy variable value that is one higher in order.

### 2.2.1 Example Problem

The executive of a company with headquarters in Saint Louis hires a chartered flight to visit all sites where the company has a location. Figure 3 shows the cities that the executive must visit. The goal is to visit all cities once and return to the headquarters in a manner that minimizes the tour cost. Hence, the chartered flight company may not necessarily visit the cities in the order of their node identifier but needs to determine the minimum distance tour.

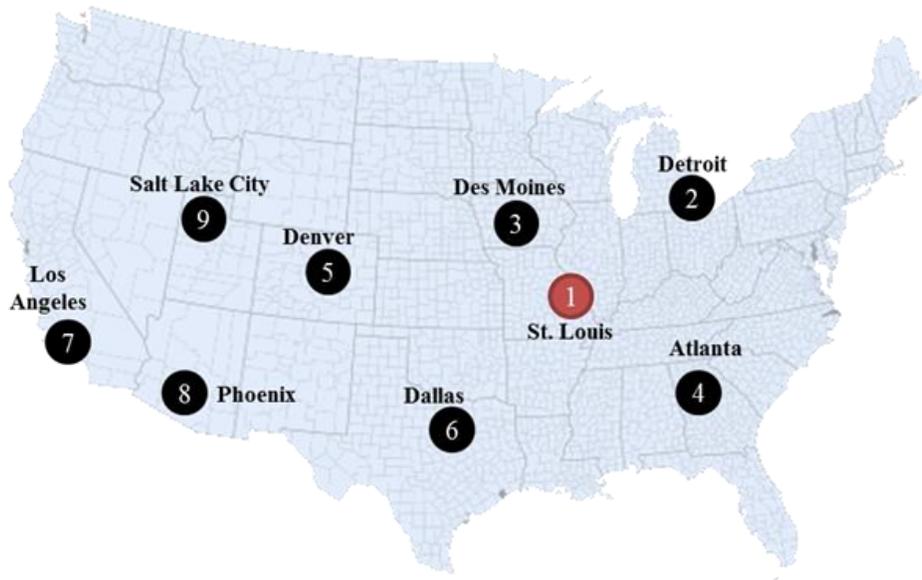

Figure 3: Minimum tour example problem.

### 2.2.2 Solution Exercise

Use GIS to generate a distance matrix for the set of cities shown in Figure 3. Set the output type to linear N × 3 format and save it into a CSV file. The course associated with this tutorial provides step-by-step instructions on how to complete this exercise. Figure 4 illustrates the solution. The values of the $u_i$ variables associated with each node is as shown. It is evident that the optimization assigned values that are equal to the position of the node in the tour. The minimum distance was approximately 4,852 miles.

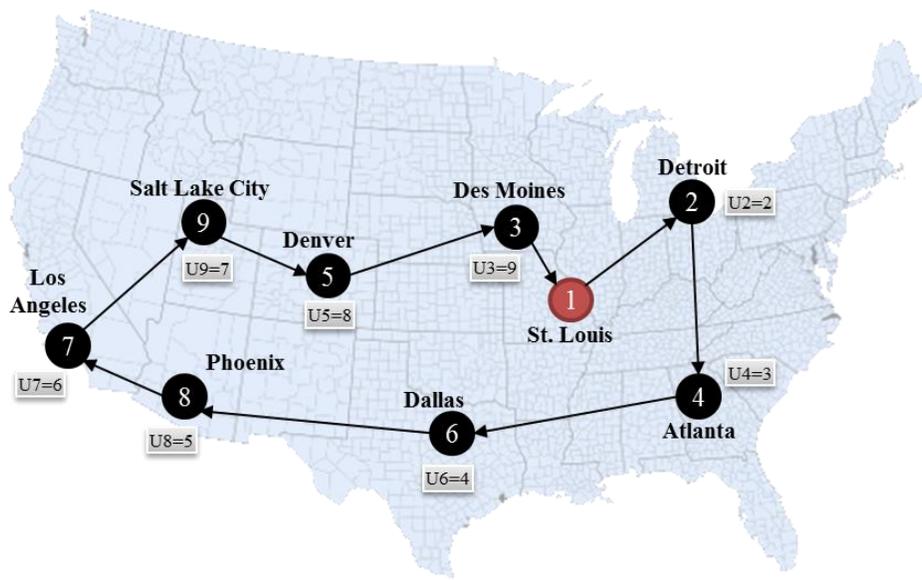

Figure 4: Solution to the minimum tour example problem.

## 2.2.3 Sensitivity Assessment

Sometimes it is not possible to complete the tour by following the minimum distance path. Figure 5 shows a scenario where the executive, while in Denver, must visit an important customer in Phoenix next. Completing the tour after that diversion resulted in a total distance of 5566.5 miles. Hence, as a sensitivity assessment of the minimum tour solution in this example, the diversion scenario resulted in nearly 15% increase in cost.

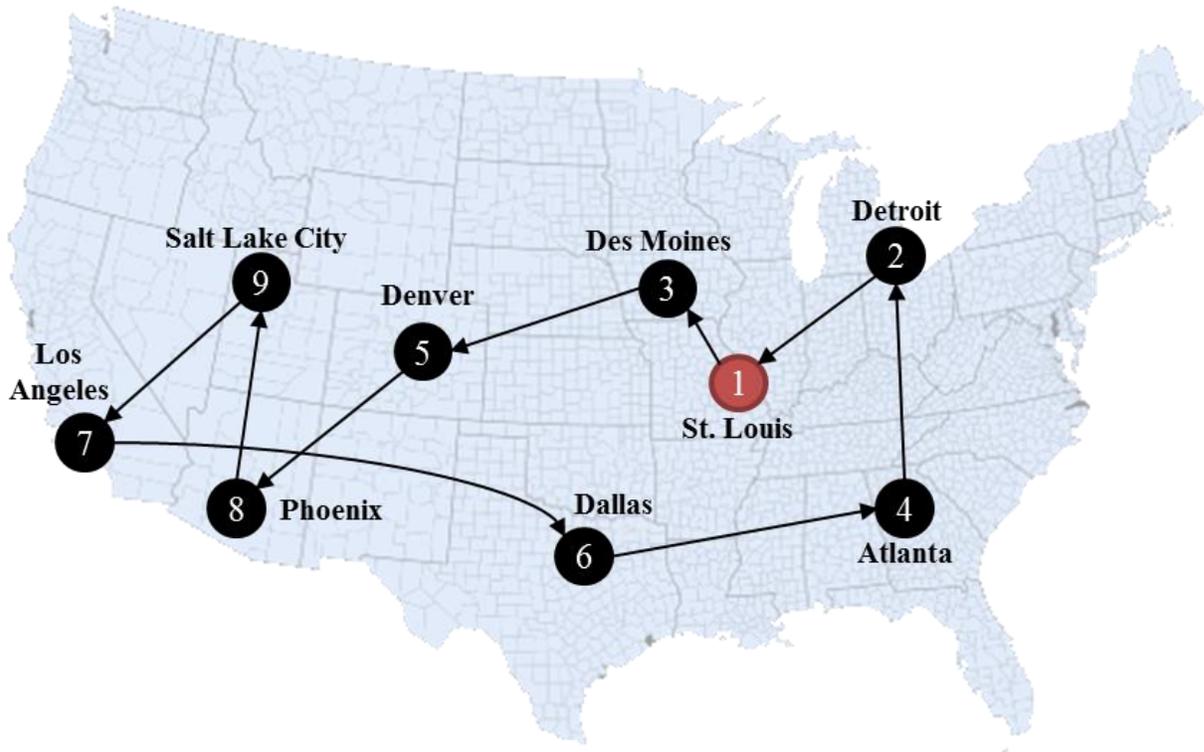

Figure 5: Sensitivity assessment of the minimum tour example problem.

## 2.2.4 Program Display

```python
# Author: Dr. Raj Bridgelall (raj.bridgelall@ndsu.edu)
# Minimum Tour Optimization
from IPython import get_ipython
get_ipython().magic('clear')                                  # Clear the console
get_ipython().run_line_magic('matplotlib', 'inline')          # plot in the iPython console
import pandas as pd
from pathlib import Path
from mip import *         # install library from Anaconda prompt: pip install mip
import numpy as np
import re
datapath_in = 'C:/Users/Admin/Documents/Minimum Cost Tour/Lab/'
infile = 'GDistance Matrix Nine Miles.csv'         # Input filename
filepath_in = Path(datapath_in + infile)           # Path name completion
df = pd.DataFrame(pd.read_csv(filepath_in, skiprows=0, usecols = (0, 1, 3))) # Read CSV to df
df = df.rename(columns={'InputID' : 'X[j]', 'TargetID' : 'X[i]'})     # Rename Columns
N = int(round(np.sqrt(df.shape[0]),0) + 1)     # Number of nodes to visit
N_Arcs = df.shape[0]                           # Possible arcs (distance matrix entries: N*(N-1))
#%% Create Variable Yij and add to model
VarNameX = []                                             # Initialize x[ij] variable names
VarU = []                                                 # Initialize u[i] variable names
for j in range(1, N+1):
    VarU.append('U'+'_'+str(j))                           # List of dummy variable names
    for i in range(1, N+1):
        if i != j:
            VarNameX.append('X'+'_'+str(i)+'_'+str(j))    # List of decision variable names
df['Xij'] = VarNameX                                      # Add variable name to table
m = Model(solver_name=CBC)                                # Instantiate optimizer
# labels and types in model vector x
x = [m.add_var(name=VarNameX[p], var_type=BINARY) for p in range(N_Arcs)]
# labels and types in model vector u
u = [m.add_var(name=VarU[q], var_type=INTEGER, lb = 2, ub = N) for q in range(1, N)]
#%% Add the cost parameters and the objective
c_parms = df.iloc[:,-2]                                   # Extract distance parameters from dataframe
m.objective = minimize(xsum(c_parms[p]*x[p] for p in range(N_Arcs)))    # Objective function
#%% Single-Entry Constraint (for each j node entered, N-1 i nodes to enter from)
# Each block of j has N-1 i's AND there are N blocks of j's
for j in range(N):
    m.add_constr( xsum(x[j*(N-1) + i] for i in range(N-1)) == 1 )
#%% Single-Exit Constraint (for each i node exited, N-1 j nodes to enter)
for i in range(1, N+1):
    # Scan all variables and extract indices of all j's associated with i
    x_idx = [k for k, s in enumerate(VarNameX) if str(i) in re.split('_',s)[1] ]
    m.add_constr( xsum(x[ x_idx[j] ] for j in range(N-1)) == 1 )
#%% Sub-tour elimination constraints
for p in range(N_Arcs):
    idx_i = int(re.split('_',VarNameX[p])[1])    # Get the i index of all decision variables
    idx_j = int(re.split('_',VarNameX[p])[2])    # Get the j index of all decision variables
    if (idx_i != 1 and idx_j != 1):
        m.add_constr( u[idx_i - 2] - u[idx_j - 2] + x[p] * N <= N - 1 )  # u[0] = "U2"
Status = m.optimize()
#%% Print the results
print('Model has {} vars, {} constraints and {} nzs'.format(m.num_cols, m.num_rows, m.num_nz))
ArcName = [VarNameX[p] for p in range(N_Arcs) if x[p].x != 0]
ArcVal = [x[p].x for p in range(N_Arcs) if x[p].x != 0]
print('Tour Arcs: {} '.format(ArcName))
print('Arc Vals: {} '.format(ArcVal))
print('Minimum Total Distance = {}'.format(m.objective_value))
print("Number of Solutions = ", m.num_solutions)
print("Status = ", Status)
df['Xij_x'] = [ x[p].x for p in range(N_Arcs) ] # Add Yij and solution to the data table
print('Confirm Sum of Distance = Objective = ', sum(df[df.Xij_x != 0].Miles)) # Total distance
for non-zero assignments
outfile = 'Minimum Tour.csv'                         # Table: demand sites covered by each server
filepath_out = Path(datapath_in + outfile)           # Full path name
df.to_csv(filepath_out, index = True, header = True) # Write CSV with the index column and header
```

Some portions of the code are compact and may benefit from some additional explanation. For example, the following code snippet establishes the single exit constraint:

```
#%% Single-Exit Constraint (for each i node exited, N-1 j nodes to enter)
for i in range(1, N+1):
    # Scan all variables and extract indices of all j's associated with i
    x_idx = [k for k, s in enumerate(VarNameX) if str(i) in re.split('_',s)[1] ]
    m.add_constr( xsum(x[ x_idx[j] ] for j in range(N-1)) == 1 )
```

The FOR loop iterates from nodes *i* ranging in index from 1 to *N*. The "N+1" is a Python standard that the range function returns a list of ordered values that is one less than the specified end index b as in range(a, b). This makes for compact syntax such as when the first index starts with 0, range(0, b) or simply range(b) returns a set of b indices starting with 0 and ending with b-1. The optimization model stores variables in an order that starts with index 0. Hence, enumerate(VarNameX) returns an iterative array containing the strings s of the variable names and the index k in the *model* array. The logic test checks to see if the second character after the "_" string separator is the same as the variable index *i*. If so, it returns the position index k in the *model* array and builds a list of those indices as x_idx using the "list compression" feature of Python. Consequently, the final list contains the indices in the *model* of all the *j* decision variables associated with that *i* instance. Subsequently, constructing the exit constraint for each node *i* references all the *j* nodes by their decision variable *x*[i, j] at the position *index* stored in the model as x[*index*]. The regular expression split function from the re library produces an array of strings that the separator character splits. In this case, there are three strings with indices from {0, 1, 2}. The string with index 1 is the *i* value.

### 2.2.5 Further Reading

1. Pataki, Gábor. "Teaching integer programming formulations using the traveling salesman problem." *SIAM review* 45, no. 1 (2003): 116-123 (Pataki, 2003).

# Chapter 3 - Spatial Optimization

Typical problems in planning and logistics that involve spatial optimization include neighborhood coverage, flow capturing to measure network operations, zone heterogeneity to avoid conflicts or competition, and to provide full-service coverage at minimum cost. The next sections discuss each type of problem.

## 3.1 Neighborhood Coverage

The general problem in this category is *set covering* where the goal is to identify a minimum set of neighborhoods to place something so that the selected locations will be adjacent to all neighboring areas. Examples include the placement of facilities such as fire stations, warehouses, and services such as emergency response centers that can serve adjacent neighborhoods within a specified amount of time. In transport logistics, sets can be neighborhoods, counties, or other areas that have boundaries. The decision variables are $X_i$ for each area $i$ covered and they take on binary values. The optimization problem is:

Minimize

$$U = \sum_{i=1}^{N} c_i X_i \tag{29}$$

Subject to:

$$\sum_{i \in A_j} X_i \geq 1 \quad (j = 1, 2, \ldots, N) \tag{30}$$

where

$$x_i \in \{0,1\} \tag{31}$$

The decision makers can normalize the cost to place a unit in an area by assigning a relative cost to the coefficient $c_i$. Table 5 is the PFT to solve the set covering problem.

Table 5: General form PFT for the Area Coverage Problem

| Area | Area *j* Covers (Constraint) | | | | Relative Cost |
|---|---|---|---|---|---|
| | 1 | 2 | … | N | |
| $X_1$ | $\alpha_{11}$ | $\alpha_{12}$ | | $\alpha_{1K}$ | 1 |
| $X_2$ | $\alpha_{21}$ | $\alpha_{22}$ | | $\alpha_{2K}$ | 1 |
| ⋮ | ⋮ | ⋮ | ⋮ | ⋮ | ⋮ |
| $X_N$ | $\alpha_{N1}$ | $\alpha_{N2}$ | | $\alpha_{NK}$ | 1 |
| ≥ | 1 | 1 | | 1 | U |

The number of areas in the "Area *j* Covers" constraint columns is equal to the number of areas or decision variables in the rows. Each column $A_j$ of the constraint matrix $A$ defines the set of areas $i$ that area *j* shares a boundary with. That is, for each constraint column the parameter is 1 if area *j* covers area *i* and zero otherwise. That is, constraint column *j* represents the set of areas covered by a placement in area *j*. For example, if area 1 covers itself plus areas 2 and 3, then the first three rows of column 1 equals 1 and the remaining rows of the column equals 0. The constraint for each column is an inequality greater than or equal to 1. Hence, the dot product of the decision variables with column *j* must be *at least* one. That is, the constraint requires at least one placement to cover an area. The cost coefficients may be absolute cost units or relative to the minimum cost of a given neighborhood so that the objective function can make placements that minimize the total cost $U$ to cover all areas with at least one placement.

### 3.1.1 Example Problem

A city with neighborhood boundaries shown in Figure 6 needs to build emergency response centers that can cover all neighborhoods. Decision makers determined that any station built in one neighborhood can serve neighborhoods that share its boundary. In this scenario, the constraint must consider that the cost to build a station in the waterfront areas of 1, 4, and 7 is twice that to build elsewhere.

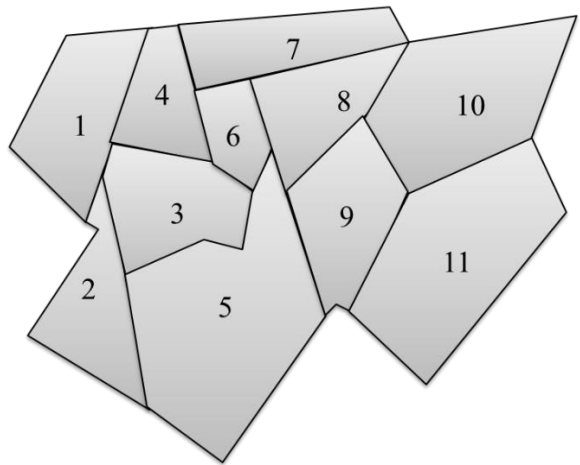

Figure 6: City neighborhood boundary for the example problem.

Table 6 shows the PFT for the problem.

Table 6: PFT for the Area Coverage Problem Exercise

| Area | Area *j* Covers | | | | | | | | | | | Cover-Items (Relative Cost) |
|---|---|---|---|---|---|---|---|---|---|---|---|---|
| | 1 | 2 | 3 | 4 | 5 | 6 | 7 | 8 | 9 | 10 | 11 | |
| $X_1$ | 1 | 1 | 1 | 1 | | | | | | | | **2** |
| $X_2$ | 1 | 1 | 1 | | 1 | | | | | | | 1 |
| $X_3$ | 1 | 1 | 1 | 1 | 1 | 1 | | | | | | 1 |
| $X_4$ | 1 | | 1 | 1 | | 1 | 1 | | | | | **2** |
| $X_5$ | | 1 | 1 | | 1 | 1 | | 1 | 1 | | | 1 |
| $X_6$ | | | 1 | 1 | 1 | 1 | 1 | 1 | | | | 1 |
| $X_7$ | | | | 1 | | 1 | 1 | 1 | | | | **2** |
| $X_8$ | | | | | 1 | 1 | 1 | 1 | 1 | 1 | | 1 |
| $X_9$ | | | | | | 1 | | 1 | 1 | 1 | 1 | 1 |
| $X_{10}$ | | | | | | | | 1 | 1 | 1 | 1 | 1 |
| $X_{11}$ | | | | | | | | | 1 | 1 | 1 | 1 |
| $\geq$ | 1 | 1 | 1 | 1 | 1 | 1 | 1 | 1 | 1 | 1 | 1 | $U$ |

The *A* matrix defines the structure of the map by indicating the bordering areas for each area with a value equal to 1. Empty values are zeros.

### 3.1.2 Solution Exercise

Figure 7a shows the solution for the problem where the cost is the same in any area. Figure 7b shows the solution when the relative cost for a waterfront area is twice the cost of any other area. The minimum number of locations is three in either case.

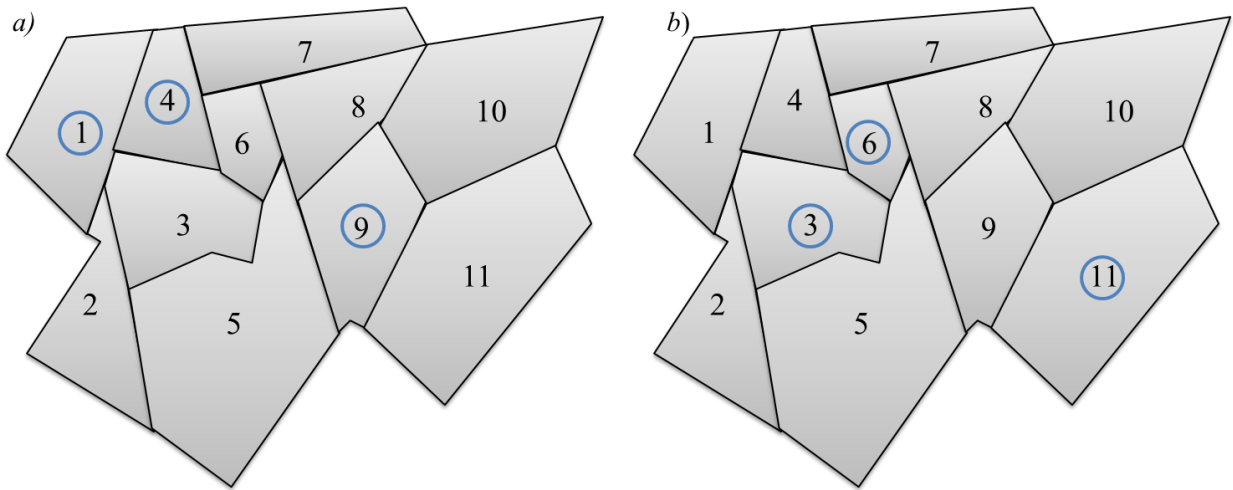

Figure 7: Solution a) for identical cost b) for higher waterfront cost.

### 3.1.3 Program Display

```
# Author: Dr. Raj Bridgelall (raj.bridgelall@ndsu.edu)
# Set covering
from IPython import get_ipython
get_ipython().magic('clear')                          # Clear the console
get_ipython().run_line_magic('matplotlib', 'inline')  # plot in the iPython console
import pandas as pd
from pathlib import Path
from mip import *      # install library from Anaconda prompt: pip install mip
import numpy as np
#%%
datapath_in = 'C:/Users/Admin/Documents/Spatial Coverage Optimization/Lab/'
infile = 'Area Coverage PFT.xlsx'                     # Input filename
filepath_in = Path(datapath_in + infile)              # Path name for untruncated signal
df = pd.DataFrame(pd.read_excel(filepath_in, skiprows=1)) # Read Problem Formulation Table (PFT)
#%%
Nc = df.shape[1]-2     # Number of resource constraints
Nd = df.shape[0]-1     # Number of decision variables
#%%
df = df.fillna(0)                                     # replace all NaN (blanks in Excel) with zeros
VarName = df.iloc[0:Nd,0]                             # String of variable names
A_parms = df.iloc[0:Nd,1:Nc+1].astype(int)            # Extract the constraints cols
b_parms = list(df.values[-1][1:-1].astype(int))       # Get constraint parameters in an array
c_parms = df.values.transpose()[-1][:-1].astype(int)  # Get objective parameters in an array
#c_parms = [1 for i in range(11)]                     # Equal cost scenario
#%%
m = Model(solver_name=CBC)                            # use GRB for Gurobi
# define list of decision vars plus store in x
x = [m.add_var(name=VarName[i],var_type=BINARY) for i in range(Nd)]
m.objective = minimize(xsum(c_parms[i]*x[i] for i in range(Nd))) # add objective function
#%% Add each constraint column
for j in range(Nc):
    m.add_constr( xsum(A_parms.iloc[i,j]*x[i] for i in range(Nd)) >= b_parms[j], "Cons"+str(j) )
```

```
#%%
Status = m.optimize()
#%% Print the results
print('Model has {} vars, {} constraints and {} nzs'.format(m.num_cols, m.num_rows, m.num_nz))
selected = [VarName[i] for i in range(Nd) if x[i].x != 0]
print('Selected Areas: {} '.format(selected))
print('Minimum Cost = {} Areas'.format(m.objective_value))
for j in range(Nc):
    print(m.constrs[j])
print("Number of Solutions = ", m.num_solutions)
print("Status = ", Status)
```

## 3.2 Flow Capturing

This optimization problem determines the optimum set of nodes for a finite number of *placements* that will capture the maximum flow through a network. Applications in transportation include *placement* of sensors to measure traffic volume, weigh-in-motion (WIM) stations to measure truck weights, driving-under-the influence (DUI) checkpoints, gasoline stations, electric vehicle charging stations, signage, and more. The problem statement is:

Maximize

$$F = \sum_{r \in R}^{N} f_r Y_r \tag{32}$$

Subject to:

$$Y_r \leq \sum_{j \in r} X_j \quad \forall r \in \mathbf{R} \tag{33}$$

and

$$\sum_{j=1}^{N} X_j = p \tag{34}$$

where

$$Y_r = \begin{cases} 1 & \text{if at least one placement is located on path } r \\ 0 & \text{otherwise} \end{cases} \tag{35}$$

$$X_j = \begin{cases} 1 & \text{if a placement is located at node } j \\ 0 & \text{otherwise} \end{cases} \tag{36}$$

The variable *r* is an index into the set of all paths **R** through the network from a start node to a terminal node. The variable $f_r$ is the *maximum* flow on path *r* of the network. A path consists of a set of arcs in the network that is a pathway from the starting node to the terminal node. Hence, some paths may have arcs in common. In the network diagram of Figure 8, flows on arcs shown by dashed arrows are not part of the analysis, but they explain why flows on some arcs can be greater than the sum of flows from the arcs that connect to them. The number of available placements is *p*.

The first set of constraints assure that the path variable $Y_r$ will be set to 1 if the optimization makes a placement on any node $X_j$ along that path. To enter the data into a PFT, the constraint must be in the standard form, which is achieved by bringing the summation to the left side of the inequality such that

$$Y_r - \sum_{j \in r} X_j \leq 0 \quad \forall r \in \mathbf{R} \tag{37}$$

The second constraint sets the available placements to the value $p$.

### 3.2.1 Example Problem

A consulting company is bidding for a grant to evaluate the performance of a new type of WIM device for the national highway system (NHS) and wants to minimize their equipment and installation budget to be cost competitive in their proposal. Figure 8 shows the network on the NHS that the proposal solicitation specified. The arc labels are the *normalized* daily truck volumes. The proposal restricts placement of a WIM station at the starting and terminal nodes.

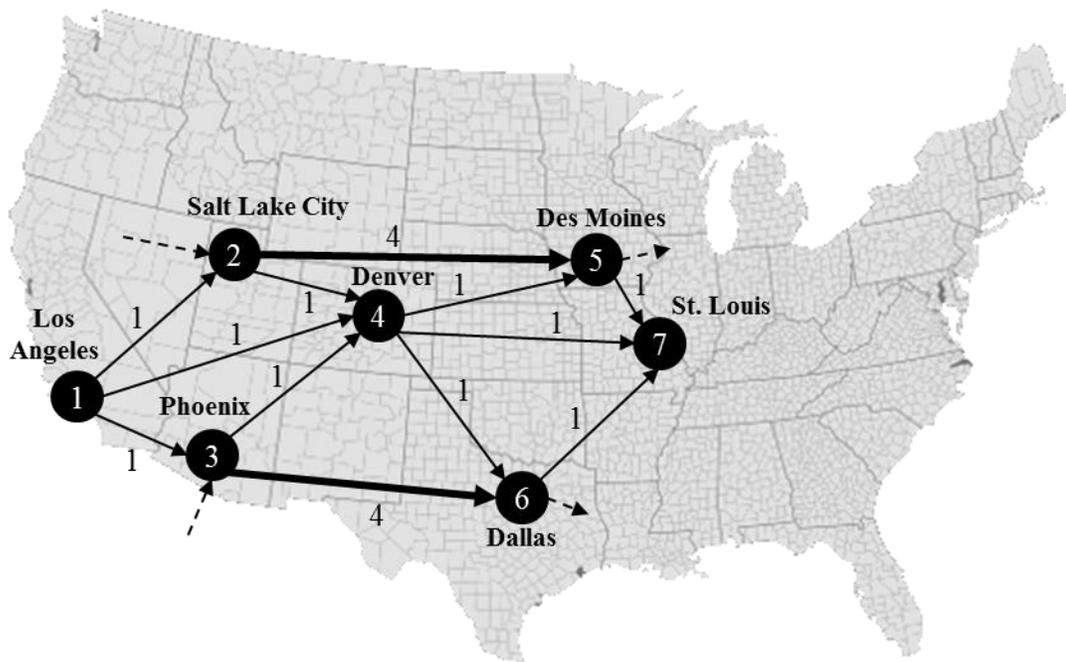

Figure 8: Flow network for the example problem.

### 3.2.2 Solution Exercise

Table 7 summarizes all paths between the start and terminal nodes. Table 8 is the corresponding PFT. The program will delete rows containing the start and terminal nodes because those nodes cannot be candidates for a WIM station. The first constraint labeled "$p$" in the header assures that the number of WIM stations sum to the total $p$ proposed. The nine path constraints $Y_r$ encode the possible paths that exist between the start and terminal nodes, which encodes the first set of constraints. The "Maximum Flow" column sets the maximum flow along each path $Y_r$ to create the objective function. The node variables are not present in the objective function, so the flow parameters are set to zero to exclude them.

Table 7: All Possible Routes from Start to Terminal Node

| Route | Path | Maximum Flow |
|---|---|---|
| Y1 | 1-2-5-7 | 4 |
| Y2 | 1-2-4-5-7 | 1 |
| Y3 | 1-4-5-7 | 1 |
| Y4 | 1-4-7 | 1 |
| Y5 | 1-4-6-7 | 1 |
| Y6 | 1-3-4-5-7 | 1 |
| Y7 | 1-3-4-7 | 1 |
| Y8 | 1-3-4-6-7 | 1 |
| Y9 | 1-3-6-7 | 4 |

Table 8: PFT for the Flow Capturing Problem

| Nodes/ Paths | \multicolumn{10}{c}{Path Constraints} | | | | | | | | | | Maximum Flow |
|---|---|---|---|---|---|---|---|---|---|---|---|
| | p | 1 | 2 | 3 | 4 | 5 | 6 | 7 | 8 | 9 | |
| $X_1$ | 1 | -1 | -1 | -1 | -1 | -1 | -1 | -1 | -1 | -1 | 0 |
| $X_2$ | 1 | -1 | -1 | | | | | | | | 0 |
| $X_3$ | 1 | | | | | | -1 | -1 | -1 | -1 | 0 |
| $X_4$ | 1 | | -1 | -1 | -1 | -1 | -1 | -1 | -1 | | 0 |
| $X_5$ | 1 | -1 | -1 | -1 | | | -1 | | | | 0 |
| $X_6$ | 1 | | | | | -1 | | | -1 | -1 | 0 |
| $X_7$ | 1 | -1 | -1 | -1 | -1 | -1 | -1 | -1 | -1 | -1 | 0 |
| $Y_1$ | | 1 | | | | | | | | | 4 |
| $Y_2$ | | | 1 | | | | | | | | 1 |
| $Y_3$ | | | | 1 | | | | | | | 1 |
| $Y_4$ | | | | | 1 | | | | | | 1 |
| $Y_5$ | | | | | | 1 | | | | | 1 |
| $Y_6$ | | | | | | | 1 | | | | 1 |
| $Y_7$ | | | | | | | | 1 | | | 1 |
| $Y_8$ | | | | | | | | | 1 | | 1 |
| $Y_9$ | | | | | | | | | | 1 | 4 |
| =, ≤ | 3 | 0 | 0 | 0 | 0 | 0 | 0 | 0 | 0 | 0 | F |

Run the program to get the solution shown in Figure 9. The solution is to place WIM stations at nodes 2, 3, and 4 to cover all routes. The maximum flow is 15 units as verified by the sum of flows for all paths $Y_r$ that those nodes contain. Change the number of units to 4 and verify that the maximum flow is still 15. Change then number of units to 2 and verify that the optimization cannot achieve node placements to a cover the maximum flow of 15.

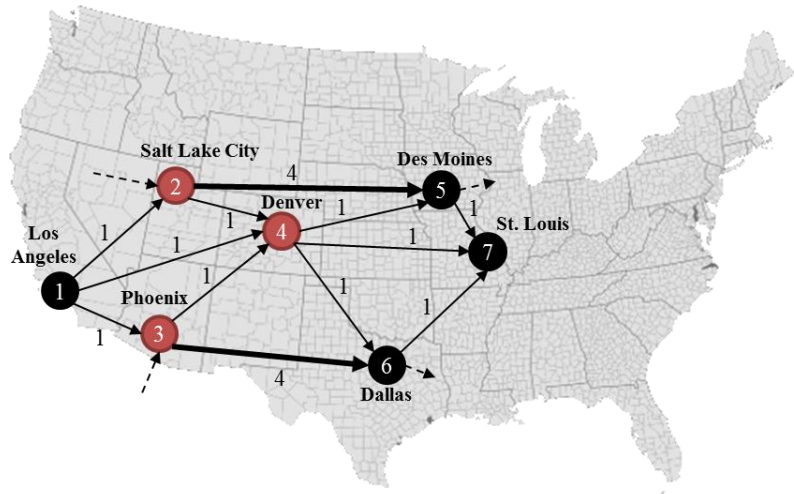

Figure 9: Flow network solution.

### 3.2.3 Program Display

```
# Author: Dr. Raj Bridgelall (raj.bridgelall@ndsu.edu)
# Flow Capture Optimization
from IPython import get_ipython
get_ipython().magic('clear')                            # Clear the console
get_ipython().run_line_magic('matplotlib', 'inline')    # plot in the iPython console
import pandas as pd
from pathlib import Path
from mip import *        # install library from Anaconda prompt: pip install mip
import numpy as np
#%%
datapath_in = 'C:/Users/Admin/Documents/Flow Capturing/Lab/'
infile = 'Flow Capture.xlsx'                            # Input filename
filepath_in = Path(datapath_in + infile)                # Path name for untruncated signal
df = pd.DataFrame(pd.read_excel(filepath_in, skiprows=1)) # Read Problem Formulation Table (PFT)
df = df.fillna(0)                                       # replace all NaN (blanks in Excel) with zeros
#%%
Nodes = int(df.values[df.shape[0]-1][0])                # Get the number of inner node variables
Facilities = int(df.values[df.shape[0]-1][1])           # Get number of facilities to deploy
df = df.drop([0, Nodes+1]).reset_index(drop=True)       # Delete start and end nodes, reset table index
#%%
Nc = df.shape[1]-2      # Number of constraints
Nd = df.shape[0]-1      # Number of decision variables
#%%
VarName = df.iloc[0:Nd,0]                               # String of variable names
A_parms = df.iloc[0:Nd,1:Nc+1].astype(int)              # Extract the constraints cols
c_parms = df.values.transpose()[-1][:-1].astype(int)    # Get objective parameters in an array
#%%
m = Model(solver_name=CBC)                              # Instantiate optimizer
x = [m.add_var(name=VarName[i],var_type=BINARY) for i in range(Nd)] # define list of decision
vars plus store in x
m.objective = maximize(xsum(c_parms[i]*x[i] for i in range(Nd)))    # add objective function
#%% Add each constraint column
m.add_constr( xsum(A_parms.iloc[i,0]*x[i] for i in range(Nd)) == Facilities, "Facilities" )
for j in range(1, Nc):
    m.add_constr( xsum(A_parms.iloc[i,j]*x[i] for i in range(Nd)) <= 0, "Cons"+str(j) )
#%%
Status = m.optimize()
#%% Print the results
print('Model has {} vars, {} constraints and {} nzs'.format(m.num_cols, m.num_rows, m.num_nz))
Facility_Loc = [VarName[i] for i in range(Nodes) if x[i].x != 0]
Covered_Routes = [VarName[i] for i in range(Nodes, Nd) if x[i].x != 0]
print('Fraction of Routes Covered: {} '.format(round(len(Covered_Routes)/(Nc-1),2)))
print('Facility Locations: {} '.format(Facility_Loc))
print('Covered Routes: {} '.format(Covered_Routes))
```

```
print('Maximum Flow = {}'.format(m.objective_value))
for j in range(Nc):
    print(m.constrs[j])
print("Number of Solutions = ", m.num_solutions)
print("Status = ", Status)
```

### *3.2.4 Further Reading*

1. Teodorović, Dušan, Milica Šelmić, Manju V. Saraswathy, Kuncheria P. Isaac, Dušan Fister, Janez Kramberger, and A. K. Sarkar. "Locating flow-capturing facilities in transportation networks: a fuzzy sets theory approach." *International Journal for Traffic & Transport Engineering* 3, no. 2 (2013): 103-111. (Teodorovic & Selmic, 2013)

## 3.3 Zone Heterogeneity

Location heterogeneity means that certain aspects of neighboring locations are different. The zone heterogeneity problem is to place items in locations that are dislike those of neighboring locations. Common reasons are to prevent interference, competition, contamination, or to yield some type of diversification. Placements can be anything from cell towers to major retail centers. A generalization of the problem is set coloring such that no adjacent area has the same color. For applications in transportation, each color could represent a facility type, a service, or a radio frequency channel, among other things. The latter is a frequent problem that requires determining how to select among a few available channels of narrow-band long-distance wireless communications to avoid radio-frequency interference. That is, devices in neighboring cells cannot use the same channel. In this case, unique channels are associated with unique colors of the generalized optimization problem.

The zone heterogeneity problem is similar in concept to the facility locating problem. However, there is a significant difference in that bordering areas (neighbors) must have different or non-competing features. Hence, the solution for *neighborhood coverage* as shown in Figure 7 will violate the zone heterogeneity objective because of the bordering facilities of (1, 4) and (3, 6). The problem formulation for zone heterogeneity requires defining subsets *A* that contain all area *pairs* $X_{ij}$ that share a boundary. The optimizer sets the decision variable $X_{ik}$ to 1 if area *i* gets color *k*. The optimization problem is, given *K* available colors and *N* nodes, minimize the number of colors assigned as:

Minimize

$$Z = \sum_{i=1}^{N} \sum_{k=1}^{K} X_{ik} \tag{38}$$

Subject to:

$$\sum_{k=1}^{K} X_{ik} = 1 \qquad (i = 1, 2, \dots, N) \tag{39}$$

and

$$X_{ik} + X_{jk} \leq 1 \qquad i,j \in A, (k = 1, 2, \dots, K) \tag{40}$$

where

$$X_{ik} \in \{0,1\} \qquad (41)$$

The first constraint of equation (39) establishes that the optimizer must assign every area one and only one color from *K* available colors. The second constraint prevents the optimizer from assigning bordering areas with the same color *k*. The final equation is a bound specifying that the decision variables are binary.

### 3.3.1 Example Problem

This example uses the same neighborhood boundary map as the neighborhood coverage problem of Figure 10. The decision variables represent each of the 11 areas assigned a color *k*. The set of areas with common borders are $A$ = {A1_2, A1_3, A1_4, A2_3, A2_5, A3_4, A3_5, A3_6, A4_6, A4_7, A5_6, A5_8, A5_9, A6_7, A6_8, A7_8, A8_9, A8_10, A9_10, A9_11, A10_11}. Table 9 is the PFT encoding the facts of the problem. The 21 columns after the column of decision variables for color *k* represents the 21 common boundary constraints. There should be a PFT for each available color *k*. The program accommodates this with an iteration loop that is equal to the number of colors available. The user can change the number of available colors manually until the optimization converges by monitoring the output status of the optimizer. The PFT does not show the constraint of assigning one and only one color to an area—the program accounts for this.

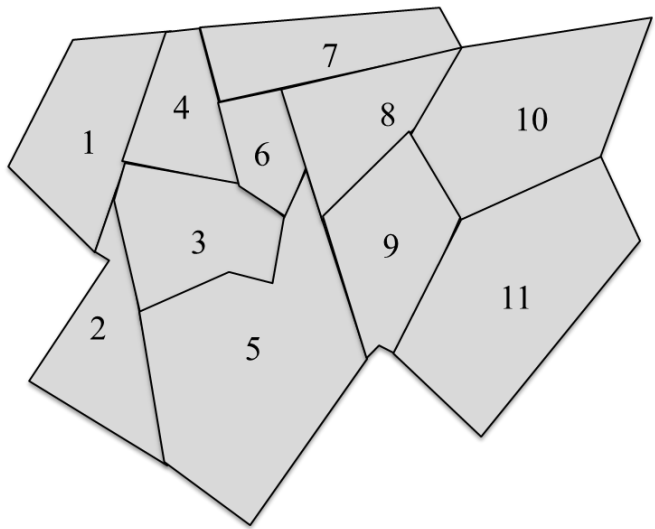

Figure 10: Neighborhood boundary map for the zone heterogeneity example problem.

Table 9: PFT for the Zone Heterogeneity Problem Exercise

| Aik | \multicolumn{19}{c|}{Pairwise Neighbor Constraints for one *k* value} | C |

| Aik | 1,2 | 1,3 | 1,4 | 2,3 | 2,5 | 3,4 | 3,5 | 3,6 | 4,6 | 4,7 | 5,6 | 5,8 | 5,9 | 6,7 | 6,8 | 7,8 | 8,9 | 8,10 | 9,10 | 9,11 | 10,11 | C |
|---|---|---|---|---|---|---|---|---|---|---|---|---|---|---|---|---|---|---|---|---|---|---|
| $X_{1k}$ | 1 | 1 | 1 | | | | | | | | | | | | | | | | | | | 1 |
| $X_{2k}$ | 1 | | | 1 | 1 | | | | | | | | | | | | | | | | | 1 |
| $X_{3k}$ | | 1 | | 1 | | 1 | 1 | 1 | | | | | | | | | | | | | | 1 |
| $X_{4k}$ | | | 1 | | | 1 | | | 1 | 1 | | | | | | | | | | | | 1 |
| $X_{5k}$ | | | | | 1 | | 1 | | | | 1 | 1 | 1 | | | | | | | | | 1 |
| $X_{6k}$ | | | | | | | | 1 | 1 | | 1 | | | 1 | 1 | | | | | | | 1 |
| $X_{7k}$ | | | | | | | | | | 1 | | | | 1 | | 1 | | | | | | 1 |
| $X_{8k}$ | | | | | | | | | | | | 1 | | | 1 | 1 | 1 | 1 | | | | 1 |
| $X_{9k}$ | | | | | | | | | | | | | 1 | | | | 1 | | 1 | 1 | | 1 |
| $X_{10k}$ | | | | | | | | | | | | | | | | | | 1 | 1 | | 1 | 1 |
| $X_{11k}$ | | | | | | | | | | | | | | | | | | | | 1 | 1 | 1 |
| ≤ | 1 | 1 | 1 | 1 | 1 | 1 | 1 | 1 | 1 | 1 | 1 | 1 | 1 | 1 | 1 | 1 | 1 | 1 | 1 | 1 | 1 | Z |

### 3.3.2 Solution Exercise

Figure 11 shows a solution that uses four colors. The colors represent categorical values and have no ordinal meaning. Users may select a suitable color scheme from the website ColorBrewer.org (Brewer & Harrower, 2020).

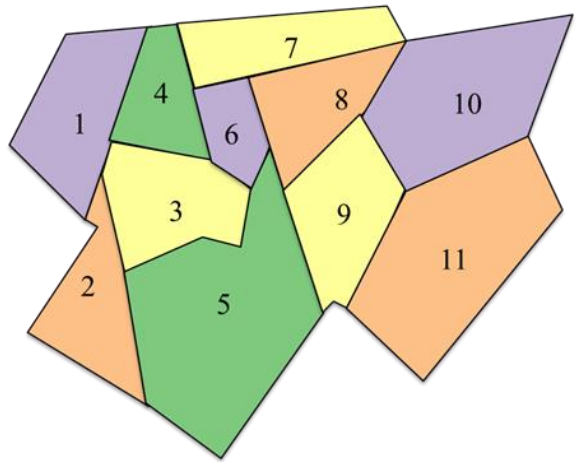

Figure 11: Solution to the zone heterogeneity example problem.

### 3.3.3 Program Display

```python
# Author: Dr. Raj Bridgelall (raj.bridgelall@ndsu.edu)
from IPython import get_ipython
get_ipython().magic('clear')                                # Clear the console
get_ipython().run_line_magic('matplotlib', 'inline')        # plot in the iPython console
import pandas as pd
from pathlib import Path
import numpy as np
from mip import *
#%%
datapath_in = 'C:/Users/Admin/Documents/Spatial Heterogeneity Optimization/Lab/'
infile = 'Color Map PFT.xlsx'                               # Input filename
filepath_in = Path(datapath_in + infile)                    # Path name for untruncated signal
df = pd.DataFrame(pd.read_excel(filepath_in, skiprows=1))   # Problem Formulation
#%%
Nc = df.shape[1]-2      # Number of boundary constraints per color
Nd = df.shape[0]-2      # Number of area decision variables per color
Ncolor = 4              # Number of colors available (modify until optimization converges)
#%%
df = df.fillna(0)                                     # replace all NaN (blanks in Excel) with zeros
#VarName = df.iloc[1:Nd+1,0]                          # String of variable names
VarName = list(df.iloc[1:Nd+1,0])                     # List of variable name strings for manipulation
A_parms = df.iloc[1:Nd+1,1:Nc+1].astype(int)          # Extract the constraints cols
#%%
m = Model(solver_name=CBC)                            # use GRB for Gurobi
#%% decision variable is whether or not to color area i with color k: create all combination
x = [m.add_var(name=VarName[i]+'_'+str(k), var_type=BINARY) for k in range(Ncolor) for i in range(Nd) ]
#%% Minimize total number of colors assigned (subject the constraints)
m.objective = minimize(xsum(x[i] for i in range(Nd*Ncolor))) # add objective function
#%% For each color, add constraints for each boundary constraint column
for k in range(Ncolor):
    for j in range(Nc):
        m.add_constr( xsum(A_parms.iloc[i,j]*x[k*Nd+i] for i in range(Nd)) <= 1,
"Cons_jk_"+str(j)+'_'+str(k) )
#%% Add constraint for a single color per area
for i in range(Nd):
    m.add_constr( xsum(x[k*Nd+i] for k in range(Ncolor) ) == 1 )
#%%
Status = m.optimize()
#%% Print the results
print('Model has {} vars, {} constraints and {} nzs'.format(m.num_cols, m.num_rows, m.num_nz))
print('Objective Function: \n',m.objective)
print('Constraints:')
for j in range(m.num_rows):
    print(m.constrs[j])
print("Status = ", Status)
print("Number of Solutions = ", m.num_solutions)
selected = [str(x[i]) for i in range(Nd*Ncolor) if x[i].x != 0]
print('Areas Colored: {} '.format(selected))
for i in range(Nd*Ncolor):
    print('{}: {} = {}'.format(i,str(x[i]),x[i].x))
#%% Number of times a color is used
C = np.zeros(Ncolor)
for k in range(Ncolor):
    for i in range(Nd):
        if x[k*Nd+i].x != 0:
            C[k] = C[k] + 1
for k in range(Ncolor):
    print('Color {} Used: {} '.format(k,C[k]))
```

The program adds constraints for each color by using the loop iterative
for k in range(Ncolor):
Indexing the decision variable as x[k*Nd+i] automatically selects the variable block for each color in the loop iteration.

## 3.4 Heterogeneity Optimization with GIS

This section demonstrates how to write a Python program that creates decision variables, constraints, and the objective function directly from the output of a GIS. Figure 12 shows the workflow.

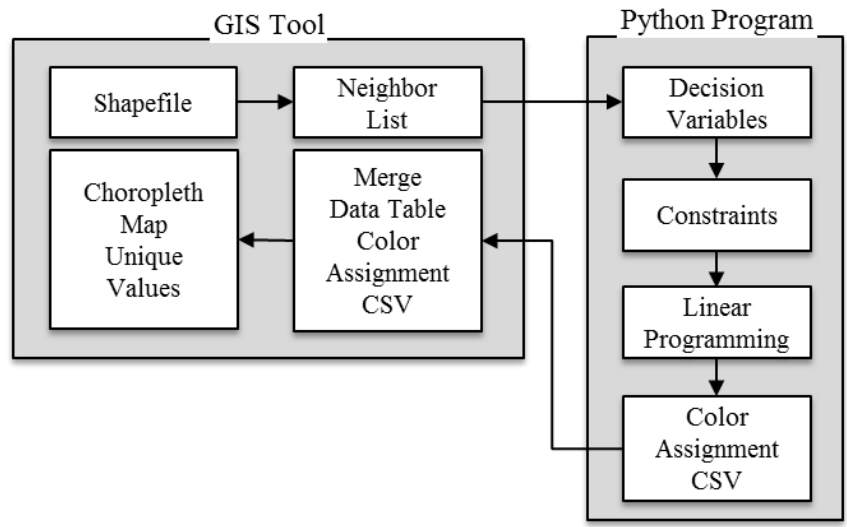

Figure 12: Workflow with GIS to solve the zone heterogeneity problem.

### 3.4.1 Exercise

The course associated with this tutorial provides a step-by-step guide to execute the workflow. The general steps are as follows:
1. Use GeoDA to import a shapefile containing all the counties of a selected state.
2. Use the GeoDA Weights Manager tool to find all the neighbors of each county and store the output in a weights file. Use the Queens weighting criteria.
3. Run the Python program in the next section to solve the spatial optimization problem and produce a CSV file of the color category assigned to each area. Change the datapath_in and infile variables to match the location and filename on your machine.
4. Merge the color assignment table with the map data table using the unique code assigned for each county.
5. Generate a categorical Choropleth map based on the unique color values.

GeoDA defines "Rook" and "Queens" contiguity to determine the criteria for neighboring polygons based on the game of chess movement rules. Neighboring polygons share only lines in the Rook criteria whereas neighbors also share vertices in the Queens criteria. Figure 13 shows the results with Texas counties using Queens contiguity weights.

### *3.4.2 Program Display*

```python
# Author: Dr. Raj Bridgelall (raj.bridgelall@ndsu.edu)
# Spatial Optimization (Color Maps) using GeoDA output
from IPython import get_ipython
get_ipython().magic('clear')                           # Clear the console
get_ipython().run_line_magic('matplotlib', 'inline')   # plot in the iPython console
import pandas as pd                                    # Data wrangling
from pathlib import Path                               # File path management
import numpy as np                                     # Numerical library
from mip import *                                      # Multiple Integer Programming Library
import re                                              # Regular Expression string processing
#%% Read in the .gal file from GeoDA
datapath_in = 'C:/Users/Admin/Documents/Zone Heterogeneity/Lab/'
infile = 'Color State_TX_Queen.gal'        # Input file from GeoDA Weights (Neighbors)
filepath_in = Path(datapath_in + infile)   # Full path name for file
list_of_lists = []                         # Initialize file input list
with open(filepath_in) as f:               # Open the file as object f
    for line in f:                                       # Get one line at a time
        inner_list = [line.strip() for line in line.split(' ')] # Split line by space dilimeter
        list_of_lists.append(inner_list)                 # Append row of string values
#%%
Ncolor = 4   # Number of available colors--increase manually from 2 until optimization converges.
N_items = int(list_of_lists[0][1])  # First line second value = totol number of areas in the map
#%% Get the decision variable names (one per item for a given color)
VarName = []                 # Initialize the list
for i in range(1,N_items*2,2): # Start at line 1, loop every other line ['Area', 'Num Neighbors']
    VarName.append(list_of_lists[i][0])     # Extract variable name string
#%% Instantiate the MIP solver
m = Model(solver_name=CBC)                                  # use GRB for Gurobi
#%% decision variable is whether or not to color area i with color k: create all combination.
# Use '_' as i_k separator. Store in object list x (optimizer refer it it only)
x = [m.add_var(name=VarName[i]+'_'+str(k), var_type=BINARY) for k in range(Ncolor) for i in range(N_items) ]
#%% Minimize total number of colors assigned (subject the constraints)
m.objective = minimize(xsum(x[i] for i in range(N_items*Ncolor)))   # Objective function
#%% For each color, add constraints for each neighbor pair in the .gal file
for k in range(Ncolor):
    for i in range(1,N_items*2,2):
        Var1 = x[VarName.index(list_of_lists[i][0]) + k*N_items]  # Get x[i] at VarName index
         for j in range(int(list_of_lists[i][1])):  # Alternate lines = ['Area', 'Num Neighbors']
            Var2 = x[VarName.index(list_of_lists[i+1][j]) + k*N_items]  # Get x[j] on next line
            m += Var1 + Var2 <= 1        # Add the x[i_k] + x[j_k] <= 1 constraint to the model
#%% Add the constraint that each area must have one and only one color
for i in range(N_items):
    m.add_constr( xsum(x[k*N_items+i] for k in range(Ncolor) ) == 1 ) # sum x[i] for all k = 1
#%% Run the optimizer and return the status
Status = m.optimize()
#%% Print the results
print('Model has {} vars, {} constraints and {} nzs'.format(m.num_cols, m.num_rows, m.num_nz))
print("Status = ", Status)
print("Number of Solutions = ", m.num_solutions)
selected = [str(x[i]) for i in range(N_items*Ncolor) if x[i].x != 0] # Variables assigned a color
print('Areas Colored: {} '.format(selected))
#%%
ItemID = [re.split('_',s) for s in selected]        # Get list of area names and color assigned
ItemDF = pd.DataFrame(ItemID)                       # Convert to data frame for CSV export
ItemDF = ItemDF.rename(columns = {0: 'County_FP', 1: 'Color'})  # Name the columns in the header
outfile = 'Color State_WY.csv'                      # Output CSV filename
filepath_out = Path(datapath_in + outfile)          # Full file name including Path
ItemDF.to_csv(filepath_out, index = False, header = True)   # Write CSV without the index column
#%% Number of times a color is used
C = np.zeros(Ncolor)                                # Initialize the color counter array
for k in range(Ncolor):
    for i in range(N_items):
        if x[k*N_items+i].x != 0:         # If a color is assigned, accumulate color count for k
            C[k] = C[k] + 1
for k in range(Ncolor):                   # Print number of time each color was used
    print('Color {} Used: {} '.format(k,C[k]))
```

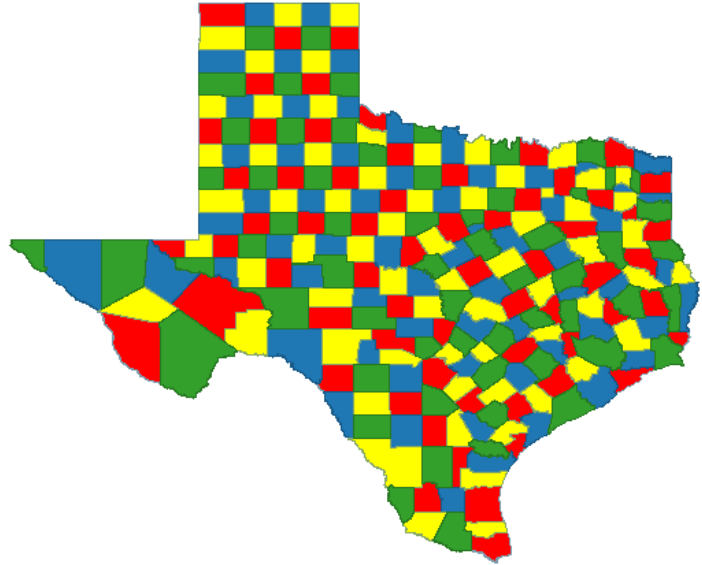

Figure 13: Results with Texas counties using Queens contiguity weights.

## 3.5 Service Coverage of Locations

This problem relates to the distribution of a finite number of *servers* such as gateways, sensors, vehicles, facilities, and services that can serve a finite set of demand locations. The problem considers where to place servers among a set of candidate locations. The same server can service multiple demand locations, but no more than one server can service the same demand location. The optimization problem is to distribute the servers in a manner that minimizes the total cost to serve all demand locations. Cost can be distance, travel cost, wireless transmission cost, and other possibilities for minimization. Using distance as a direct cost, the variables are:

*I*  the set of *N demand* node locations indexed by *i*
*J*  the set of *M* candidate server locations indexed by *j*
*p*  the number of servers available
$d_{ij}$  the physical distance between demand node *i* and candidate server location *j*.

The problem formulation is:

Minimize

$$D = \sum_{i=1}^{N} \sum_{j=1}^{M} d_{ij} Y_{ij} \qquad (42)$$

Subject to:

$$\sum_{j=1}^{M} Y_{ij} = 1, \quad \forall i \in I \qquad (43)$$

and

$$\sum_{j=1}^{M} X_j = p \tag{44}$$

and

$$Y_{ij} \leq X_j, \quad \forall i \in I, \forall j \in J \tag{45}$$

where

$$Y_{ij} = \begin{cases} 1 & \text{location } i \text{ is served from location } j \\ 0 & \text{otherwise} \end{cases}, \quad \forall i \in I, \forall j \in J \tag{46}$$

$$X_j = \begin{cases} 1 & \text{if server is placed at location } j \\ 0 & \text{otherwise} \end{cases}, \quad \forall j \in J \tag{47}$$

The objective function selects candidate sites that minimizes the overall service distance in the network. The first constraint assures that one and only one server will serve a demand site. The second constraint assures that the number of server placements is exactly $p$. The third constraint assures that if the optimizer places a server at location $j$ to serve location $i$ then it must set server $j$ location as assigned. All decision variables are binary.

### 3.5.1 Scenario Problem with QGIS

Many cities have been considering land use planning to accommodate micromobility. During the early deployment of micromobility vehicles such as electric scooters and electric bicycles, cities faced many challenges (NACTO, 2018). The convenience of stowing a vehicle anywhere resulted in clutter and interference with pedestrian traffic. Consequently, cities have been restricting the use of undocked micromobility vehicles. Meanwhile, companies have proposed new types of micromobility vehicles, including some that have an enclosed compartment to provide heat and shelter. Cities recognize that the demand for such vehicles will increase because they provide a lower-cost alternative for short trips, and they can connect to public transportation for trip completion or longer trips. Hence, affordable and accessible micromobility devices could help relieve traffic congestion, particularly if deployment results in a mode shift away from single-occupancy vehicles.

The problem scenario posed for this exercise is that one city proposed to design small, marked spaces next to bus stops and subway stations where users can safely stow micromobility vehicles after use or incur a penalty if they did not. The pilot study will deploy a fleet of micromobility vehicles near a popular park in a densely populated suburb. Planners wish to start with three bus stops selected from nine candidate locations, and the subway stations located within four kilometers of the park centroid. Their objective is to minimize the total distance between the selected locations and the subway stations. Table 10 shows the PFT for setting up the optimization problem in software.

Table 10: PFT for the Service Coverage Problem

| Service Links | p | Service Constraints for $j$ | | | | | | Distance Matrix |
|---|---|---|---|---|---|---|---|---|
| | 0 | 1 | ... | M | 1 | ... | M | |
| $Y_{11}$ | 0 | 1 | | | 1 | | | $d_{11}$ |
| ⋮ | 0 | | 1 | | | 1 | | ⋮ |
| $Y_{NM}$ | 0 | | | 1 | | | 1 | $d_{NM}$ |
| $X_1$ | 1 | | | | -1 | | | 0 |
| ⋮ | 1 | | | | | -1 | | ⋮ |
| $X_M$ | 1 | | | | | | -1 | 0 |
| $=, \leq$ | p | 1 | 1 | 1 | 0 | 0 | 0 | D |

### 3.5.2 Solution Exercise

Figure 14 illustrates the solution.

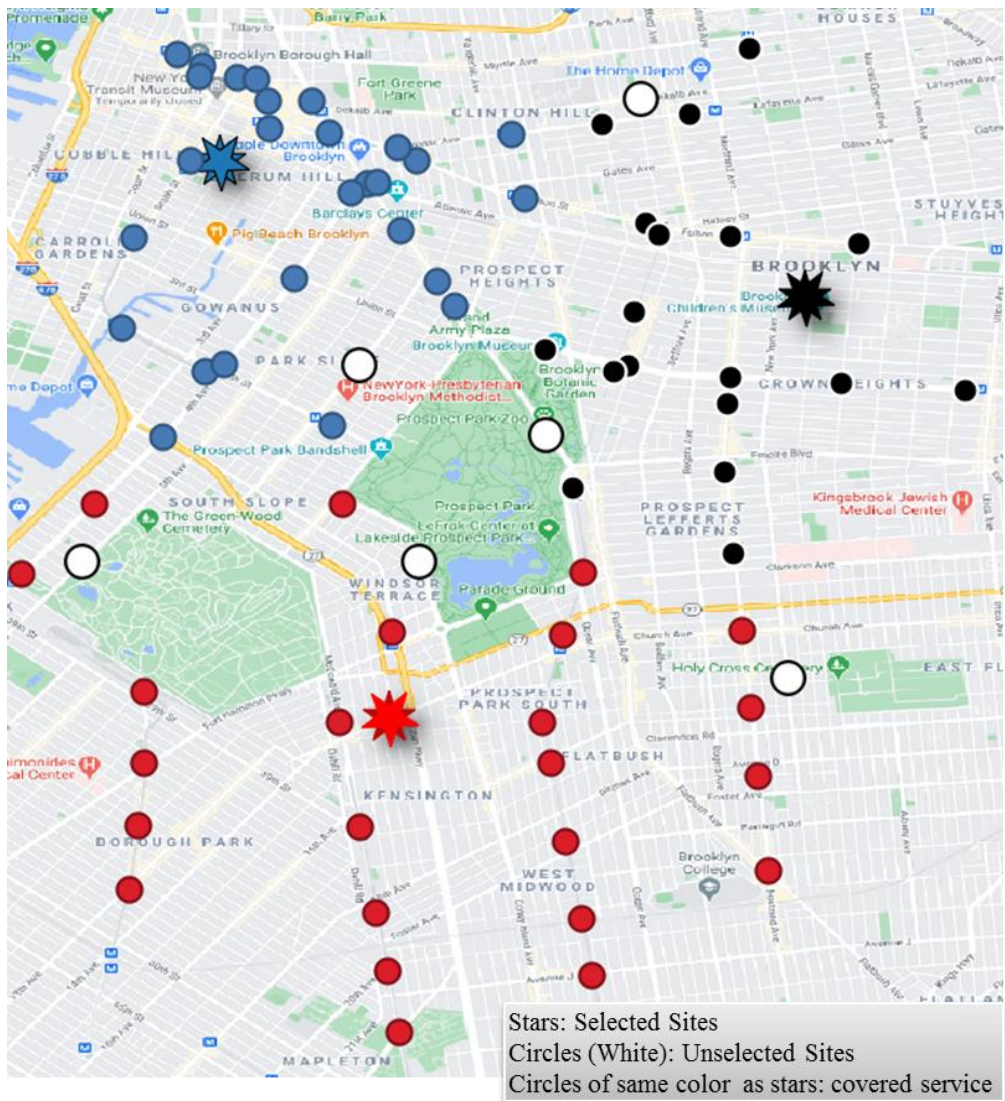

Figure 14: Solution to the service coverage example problem.

The course associated with this tutorial provides step-by-step instructions on how to complete this exercise. The general workflow is as follows:
1) Download shapefiles for New York City subway and bus stops.
2) Use QGIS to select all subway and bus stop that are within 4 kilometers from the centroid of Prospect Park in Brooklyn, New York.
3) Use QGIS to select nine candidate sites at bus stops around the park.
4) Use QGIS to produce a distance matrix for all combination of distances between the candidate sites and the subway stations.
5) Save the distance matrix as a CSV file.
6) Run the Python program to use the distance matrix as an input file, build the optimization model, and write the results into a cover map file of the subway stations covered by the selected bus-stop sites.
7) Use QGIS to display a color-coded map of the results.

### 3.5.3 Program Display

```
# Author: Dr. Raj Bridgelall (raj.bridgelall@ndsu.edu)
# Service Coverage Optimization
from IPython import get_ipython
get_ipython().magic('clear')                                    # Clear the console
get_ipython().run_line_magic('matplotlib', 'inline')            # plot in the iPython console
import pandas as pd
from pathlib import Path
from mip import *         # install library from Anaconda prompt: pip install mip
import numpy as np
#%%
datapath_in = 'C:/Users/Admin/Documents/Land Use Modeling/Lab/'
infile = 'Distance Matrix Bus Subway.csv'        # Input filename
filepath_in = Path(datapath_in + infile)         # Path name completion
df = pd.DataFrame(pd.read_csv(filepath_in, skiprows=0)) # Read CSV into dataframe
df = df.rename(columns={'InputID' : 'BusStop[j]', 'TargetID' : 'Subway[i]'})    # Rename Columns
#%%
N_Servers = 3                          # Number of serving stations
M_BusStop = 9                          # Number of target sites for serving stations
N_Subways = int(df.shape[0]/M_BusStop) # Number of demand sites to serve
MN = df.shape[0]                       # Total combination of MN for loops
#%% Create Variable Yij and add to model
VarNameY = []                                      # Initialize list
for j in range(1, M_BusStop+1):
    for i in range(1, N_Subways+1):
        VarNameY.append('Y'+str(i)+'_'+str(j))     # Create list of decision variable names
m = Model(solver_name=CBC)                         # Instantiate optimizer
y = [m.add_var(name=VarNameY[p], var_type=BINARY) for p in range(MN)] # y = y vars
#%% Add the cost parameters and the objective
c_parms = df.iloc[:,-1]                            # Extract distance parameters from dataframe
m.objective = minimize(xsum(c_parms[p]*y[p] for p in range(MN)))    # objective function
#%% constraint for demand i serviced by only one server j
for i in range(N_Subways):
    m.add_constr( xsum(y[j*N_Subways + i] for j in range(M_BusStop)) == 1 )
#%% Create Variable Xi and add to model
VarNameX = []
for i in range(1, N_Subways+1):
    VarNameX.append('X'+str(i))
x = [m.add_var(name=VarNameX[i], var_type=BINARY) for i in range(N_Subways)]   # x = x vars
#%% Add constraint Yij - Xj <= 0
for j in range(M_BusStop):
    for i in range(N_Subways):
        m.add_constr( y[j*N_Subways + i] - x[j] <= 0 )  # y has M blocks of i (each N long)
#%% Add constraint sum(xj) = p
m.add_constr( xsum( x[j] for j in range(M_BusStop) ) == N_Servers )
#%%
Status = m.optimize()
#%% Print the results
```

```python
print('Model has {} vars, {} constraints and {} nzs'.format(m.num_cols, m.num_rows, m.num_nz))
Server_Loc = [VarNameX[i] for i in range(M_BusStop) if x[i].x != 0]
Service_Net = [VarNameY[p] for p in range(MN) if y[p].x != 0]
print('Server Locations: {} '.format(Server_Loc))
print('Service Network: {} '.format(Service_Net))
print('Minimum Total Distance = {}'.format(m.objective_value))
print("Number of Solutions = ", m.num_solutions)
print("Status = ", Status)
#%% Add Yij variables and solution to the data table
df['Yij'] = [ VarNameY[p] for p in range(MN) ]
df['Yij_x'] = [ y[p].x for p in range(MN) ]
print('Confirm Sum of Distance = Objective = ', sum(df[df.Yij_x != 0].Distance)) # NZ tot dist
#%% Add Xj variables and solution to the data table
# Table has M blocks of i (each N long)
df['Xj'] = [ VarNameX[j] for j in range(M_BusStop) for i in range(N_Subways) ]
df['Xj_x'] = [ x[j].x for j in range(M_BusStop) for i in range(N_Subways) ]     # matching values
#%% Table of demand sites covered by each server
CoverMap = df.pivot_table(values='Yij_x', index = ['Subway[i]'], columns=['BusStop[j]'],
aggfunc=np.sum)
CoverTable = pd.DataFrame(CoverMap.sum().rename('Served'))      # See in 'Variable Explorer'
outfile = 'Cover Table.csv'                                     # Output filename
filepath_out = Path(datapath_in + outfile)                      # Full path name
CoverTable.to_csv(filepath_out, index = True, header = True)    # Write CSV w/ index col & header
#%% List of demand sites covered by servers
ColorMap = pd.get_dummies(CoverMap).idxmax(1)                   # Reverse one-hot-encoding
ColorMap = ColorMap.rename('BusStop[j]')                        # Rename column
outfile = 'Subway Color Map.csv'                                # Output filename
filepath_out = Path(datapath_in + outfile)                      # Full path name
ColorMap.to_csv(filepath_out, index = True, header = True)      # Write CSV w/ index col & header
```

# Chapter 4 - Spatial Logistics

Problems in spatial logics involve the distribution of items to satisfy supply and demand constraints while minimizing cost or maximizing flow. The next sections cover three types of spatial logistics problems that involve distribution, flow, and warehouse location optimization.

## 4.1 Spatial Distribution

This category of optimization problems involve decision making to meet supply and demand constraints. The decision variables $X_{ij}$ are associated with each combination of supply entity $i$ and demand entity $j$. The objective function is a *supplier* cost minimization where coefficient $c_{ij}$ represents the per-unit cost of production plus transportation. The problem formulation for a single item from $M$ suppliers and $N$ demand locations is:

Minimize

$$Z = \sum_{i=1}^{M} \sum_{j=1}^{N} c_{ij} X_{ij} \tag{48}$$

Subject to:

$$\sum_{i=1}^{M} X_{ij} = d_j \quad (j = 1, 2, \dots, N) \tag{49}$$

and

$$\sum_{j=1}^{N} X_{ij} \leq \sum_{j=1}^{N} d_j \quad (i = 1, 2, \dots, M) \tag{50}$$

where

$$X_{ij} \geq 0 \quad (i = 1, 2, \dots, M), \ (j = 1, 2, \dots, N) \tag{51}$$

The first constraint of equation (49) is a *demand* constraint assuring that the demand is met for all locations. It states that the total units supplied from all suppliers $i$ to location $j$ must be equal to the demand from location $j$. The second constraint is a *supply* constraint that establishes an upper bound on supplier $i$. It states that the total amount from supplier $i$ can be at most the total demand from all locations $j$. This constraint determines the maximum design capacity for each supplier $i$ when the logistics cost of the *system* is minimum. An alternate supply constraint is to set the existing capacity of each supplier as

$$\sum_{j=1}^{N} X_{ij} \leq U_i \quad (i = 1, 2, \dots, M) \tag{52}$$

The variable bound Equation (51) establishes a non-negative constraint on the quantity of items supplied. That is, the decision variables are *integers*, and they have only a lower bound.

### 4.1.1 Example Problem

This example is adopted from a popular example involving warehouse suppliers and stores (Mitchell, Kean, Mason, O'Sullivan, & Phillips, 2020). In this scenario, the supplier has two warehouses that supply smartphones to five stores in the region. Each week, the average demand from stores 1 through 5 are 500, 900, 1,800, 200, and 700 units, respectively. Figure 15 is a graphical representation of the problem.

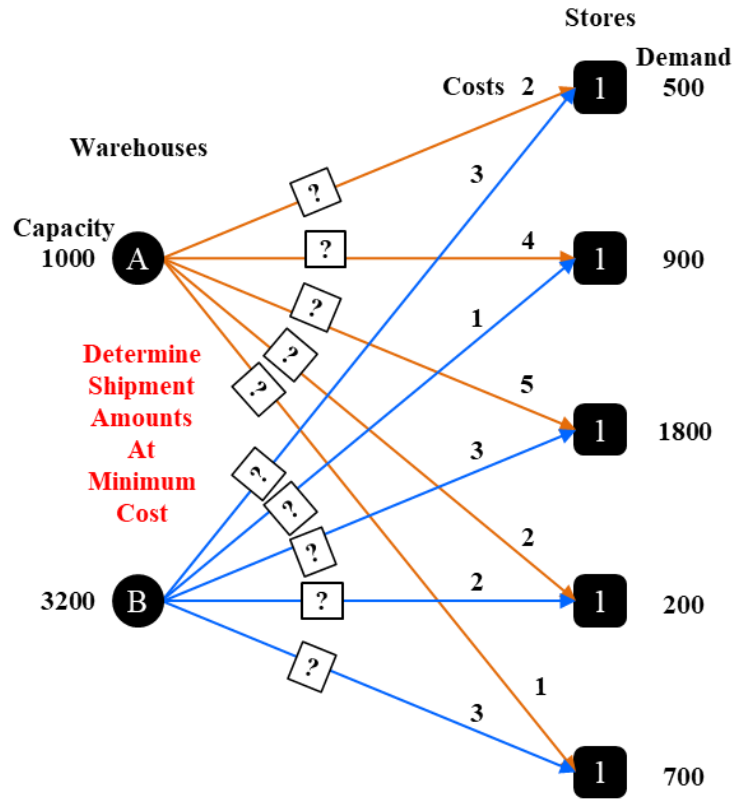

Figure 15: Graphical representation of the spatial logistics problem.

One week, the warehouse manager noted that they could supply only 1,000 units from warehouse A and 3,200 units from warehouse B. Given that the cost per unit shipment depends on the truck route between a given warehouse and a store (Table 11), the manager needed to determine the optimum distribution of shipments from each warehouse to stores that will minimize the shipping costs while meeting all demand.

Table 11: Shipping Cost for the Spatial Distribution Exercise

|  | From Warehouse | |
| --- | --- | --- |
| To Store | A ($) | B ($) |
| 1 | 2 | 3 |
| 2 | 4 | 1 |
| 3 | 5 | 3 |
| 4 | 2 | 2 |
| 5 | 1 | 3 |

Table 12 shows the PFT for this spatial distribution problem. The five demand constraints are equalities (eq) and the two supply constraints are upper bounds (ub).

Table 12: PFT for the Spatial Distribution Exercise

| Ship | Stores (Demand Constraint eq) | | | | | Warehouse (Supply Constraint ub) | | Cost |
|---|---|---|---|---|---|---|---|---|
| $i$ to $j$ | 1 | 2 | 3 | 4 | 5 | A | B | ($) |
| $X_{11}$ | 1 | | | | | 1 | | 2 |
| $X_{12}$ | | 1 | | | | 1 | | 4 |
| $X_{13}$ | | | 1 | | | 1 | | 5 |
| $X_{14}$ | | | | 1 | | 1 | | 2 |
| $X_{15}$ | | | | | 1 | 1 | | 1 |
| $X_{21}$ | 1 | | | | | | 1 | 3 |
| $X_{22}$ | | 1 | | | | | 1 | 1 |
| $X_{23}$ | | | 1 | | | | 1 | 3 |
| $X_{24}$ | | | | 1 | | | 1 | 2 |
| $X_{25}$ | | | | | 1 | | 1 | 3 |
| =, ≤ | 500 | 900 | 1800 | 200 | 700 | 1000 | 3200 | Z |

### 4.1.2 Solution Exercise

Figure 16 shows a graphical representation of the solution. The minimum shipment cost was $8,600.

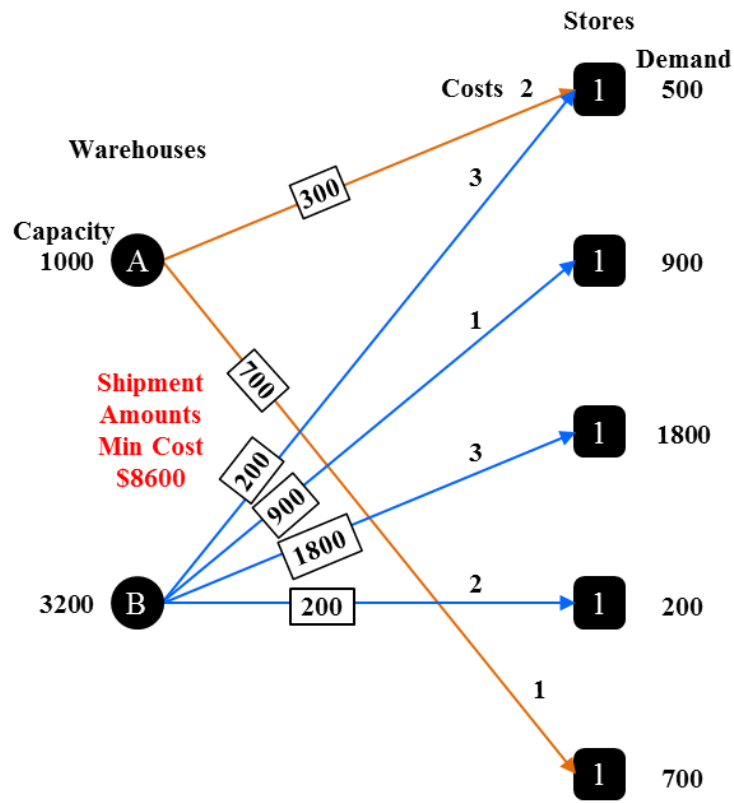

Figure 16: Representation of the spatial distribution solution.

Create a copy of the PFT in an Excel file and then run the Python program in the next subsection.

### *4.1.3 Program Display*

```python
# Author: Dr. Raj Bridgelall (raj.bridgelall@ndsu.edu)
# Spatial distribution to balance supply and demand
from IPython import get_ipython
get_ipython().magic('clear')                                  # Clear the console
get_ipython().run_line_magic('matplotlib', 'inline')          # plot in the iPython console
import pandas as pd
from pathlib import Path
from mip import *         # install library from Anaconda prompt: pip install mip
#%%
datapath_in = 'C:/Users/Admin/Documents/Spatial Logistics and Flows/Lab/'
infile = 'Warehouse Distribution PFT.xlsx'          # Input filename
filepath_in = Path(datapath_in + infile)            # Path name for untruncated signal
df = pd.DataFrame(pd.read_excel(filepath_in, skiprows=1))  # Read Problem Formulation Table (PFT)
#%%
Nc = df.shape[1]-2     # Number of constraints
Nd = df.shape[0]-1     # Number of decision variables
#%%
df = df.fillna(0)                            # replace all NaN (blanks in Excel) with zeros
VarName = df.iloc[0:Nd,0]                    # String of variable names
A_df = df.iloc[0:Nd,1:Nc+1].astype(int)      # Extract the constraints cols
A_parms_eq = A_df.iloc[:,0:Nc-2]             # Extract A_eq
A_parms_ub = A_df.iloc[:,Nc-2:]              # Extract A_ub (already in standard ub form)
b_parms = df.values[-1][1:-1].astype(int)    # Get constraint parameters in an array
b_parms_eq = b_parms[0:Nc-2]                 # Extract b_eq
b_parms_ub = b_parms[Nc-2:]                  # Extract b_ub
c_parms = df.values.transpose()[-1][:-1].astype(int)    # Get objective parameters in an array
#%%
m = Model(solver_name=CBC)                            # use GRB for Gurobi
# define list of decision vars plus store in x
x = [m.add_var(name=VarName[i],var_type=INTEGER) for i in range(Nd)]
m.objective = minimize(xsum(c_parms[i]*x[i] for i in range(Nd))) # add objective function
#%% Add each constraint column
for j in range(len(b_parms_eq)):
    m.add_constr( xsum(A_parms_eq.iloc[i,j]*x[i] for i in range(Nd)) == b_parms_eq[j],\
              "ConsD_"+str(j) )
for j in range(len(b_parms_ub)):
#    m.add_constr( xsum(A_parms_ub.iloc[i,j]*x[i] for i in range(Nd)) <= sum(b_parms_eq),\
#              "Cons_AnyCap_"+str(j) )
    m.add_constr( xsum(A_parms_ub.iloc[i,j]*x[i] for i in range(Nd)) <= b_parms_ub[j],\
              "Cons_Cap_"+str(j) )
#%%
Status = m.optimize()
#%% Print the results
print('Model has {} vars, {} constraints and {} nzs'.format(m.num_cols, m.num_rows, m.num_nz))
print("Number of Solutions = ", m.num_solutions)
print("Status = ", Status)
Arcs = [VarName[i] for i in range(Nd) if x[i].x != 0]
ArcVals = [x[i].x for i in range(Nd) if x[i].x != 0]
for i in range(len(Arcs)):
    print('{} = {} '.format(Arcs[i], ArcVals[i]))

print('Minimum Cost = ${}'.format(m.objective_value))
print('Constraints:')
for j in range(Nc):
    print(m.constrs[j])
```

As with all examples, be sure to change the variable "datapath_in" to match the path where the Excel file is located. Change the "infile" variable to match the name of the Excel file. The program extracts the PFT by skipping the first row to create a data frame. Because the software imports blank cells as "nan" it replaces those with zeros. It stores the variable names in the first column as a series of strings. Change the problem scenario to designing capacities for each warehouse to minimize the logistics cost of the system by switching the comments on the inequality constraints. Observe that the total logistics cost of the system drops to $8,400 by changing the capacities of warehouse A and B to 1,400 and 2,700, respectively.

### 4.1.4 Linear Algebra Formulation

An alternative optimizer for solving simpler problems that involve only continuous variables is available in the `scipy.optimize.linprog` library. The optimizer solves linear programming problems in the following form:

$$\min_{x} c^T x \tag{53}$$

such that

$$A_{ub} x \leq b_{ub} \tag{54}$$

$$A_{eq} x = b_{eq} \tag{55}$$

$$l \leq x \leq u \tag{56}$$

where $x$ is a vector of decision variables, $c$ is a vector of objective function coefficients, $b_{ub}$ is a vector of the *upper bound* constants on the right side of the *inequality* constraint equations, $b_{eq}$ is a vector of the equality constants on the right side of the *equality* constraint equations, $l$ and $u$ are vectors of the lower and upper bounds, respectively, of the associated decision variable in vector $x$. A program can derive the $A_{ub}$ and $A_{eq}$ matrices by transposing the corresponding $A$ matrices in the PFT. That is, each row of $A_{ub}$ and $A_{eq}$ must contain the coefficients of the linear inequality and linear equality constraints, respectively, on the decision variables in vector $x$.

Note that the standard form for the inequality constraint (54) requires the constant to be an upper bound (that is, less than or equal). Hence, the analyst must convert any constraint formulated as a lower bound (that is, greater than or equal) to an upper bound form by multiplying each side of the inequality by negative one. One limitation of this optimizer is that it has no provisions to directly accommodate binary variables. However, it is sometimes, but not always possible to simulate binary decision variables by setting the lower bound to 0 and the upper bound to 1.

### 4.1.5 Linear Algebra Formulation

The following program solves the spatial distribution problem using the alternate optimizer. It extracts the constraint parameters into a matrix and transposes it to the require form of $A_{eq}$. It also extracts the $c$ and $b_{eq}$ parameters into arrays that the optimizer requires. The program stores the bounds for each variable as a tuple with zero as the lower bound and no upper bound. The program then calls the optimizer and passes all the required parameters in their expected format. The "revised simplex" optimization method produces the best results for integer programming. The program produces the results as an object that contains the following elements:

| | |
|---|---|
| con: | an array of floating-point values of the *equality* constraint gap. |
| fun: | a floating-point value of the optimal value of the objective function. |
| message: | a string describing the exit status of the algorithm. |
| nit: | an integer of the total number of iterations performed in all phases. |
| slack: | an array of floating-point values of the *inequality* constraint gap. |

status: an integer representing a status where
0 : Optimization proceeding nominally.
1 : Iteration limit reached.
2 : Problem is infeasible.
3 : Problem is unbounded.
4 : Numerical difficulties encountered.
success: a Boolean that is true when the algorithm completes successfully.
x: an array of floating-point values containing the solution for the decision variables.

*4.1.6 Program Display*

```
# Author: Dr. Raj Bridgelall (raj.bridgelall@ndsu.edu)
# Spatial distribution using Scipy optimizer
from IPython import get_ipython
get_ipython().magic('clear')                            # Clear the console
get_ipython().run_line_magic('matplotlib', 'inline')    # plot in the iPython console
import pandas as pd
from pathlib import Path
from scipy.optimize import linprog
import numpy as np
#%%
datapath_in = 'C:/Users/Admin/Documents/Spatial Distribution/Lab/'
infile = 'Warehouse Distribution PFT.xlsx'              # Input filename
filepath_in = Path(datapath_in + infile)                # Path name for untruncated signal
df = pd.DataFrame(pd.read_excel(filepath_in, skiprows=1))   # Problem Formulation
#%%
Nc = df.shape[1]-2      # Number of resource constraints
Nd = df.shape[0]-1      # Number of decision variables
#%%
#%%
df = df.fillna(0)                                       # replace all NaN (blanks in Excel) with zeros
A_df = df.iloc[0:Nd,1:Nc+1].astype(int)                 # Extract the A parameters subset
A_parms_eq = A_df.iloc[:,0:Nc-2]                        # Extract A_eq
A_parms_eq = A_parms_eq.values.transpose()              # Convert A parm subset to matrix, transposed
b_parms = df.values[-1][1:-1].astype(int)               # Get constraint parameters in an array
b_parms_eq = b_parms[0:Nc-2]                            # Extract b_eq

A_parms_ub = A_df.iloc[:,Nc-2:]                         # Extract A_ub (already in standard ub form)
A_parms_ub = A_parms_ub.values.transpose() # Transpose A params to standard form of optimizer
b_parms_ub = b_parms[Nc-2:]                             # Extract b_ub

c_parms = df.values.transpose()[-1][:-1].astype(float)  # Get objective parameters in an array
x_bounds = (0, None)                                    # set lower and upper bounds for the variables
#%% Call the optimizer
# https://docs.scipy.org/doc/scipy/reference/generated/scipy.optimize.linprog.html
res = linprog(c_parms, A_ub=A_parms_ub, b_ub=b_parms_ub, A_eq=A_parms_eq, b_eq=b_parms_eq,
              bounds=x_bounds,
              method='revised simplex')
#%% Print the results
print(res); print("\n")
Sol = res.x.astype(int)
print("Values for Decision Variables: ", Sol)
Arcs = [df.values.transpose()[0][i] for i in range(Nd) if Sol[i] > 0]
Arc_Vals = [Sol[i] for i in range(Nd) if Sol[i] > 0]
print("Arcs: ", Arcs)
print("Arc Values: ", Arc_Vals)
print("Minimum Cost ($) = ", np.dot(Sol,c_parms))
print("Validation from LinProg = ", res.fun)
```

## 4.2 Flow Maximization

This type of problem optimizes flows through a network for a set of origins and destinations. Typical applications in transportation include:

- **Evacuation routing**: the constraints of each arc in a path can be the roadway capacity in maximum vehicle volume per hour based on the speed limit, traffic control devices, and geometric features of the roadways (number of lanes, shoulder width, medians, terrain, interchange density, etc.). The traffic source nodes can be towns affected by a threat and the sink nodes can be candidate towns that can serve as temporary sanctuaries.
- **Pipeline network throughput**: the constraints of each arc in a path can be the flow capacity in maximum flow units per day (such as barrels of oil). The source nodes can be oil production regions such as the Bakken shale, and the sink nodes can be refineries or reservoirs.
- **Railroad network throughput**: the constraints of each arc in a path can be the track capacity in maximum flow units per day (such as shipping containers). The source nodes can be shipping ports or terminals, and the sink nodes can be transshipment terminals that load trucks for delivery to warehouses or stores.

The links considered are candidate routes between nodes. The problem defines links only for paths that exist between nodes. The node labels and constraints that define the network connectivity and flow conservation of the nodes are similar to that of the shortest path problem. However, there are four major differences:
1) the objective function defines the source node $s$ by maximizing flows out of it.
2) the constraints for flow conservation define the network and the node that has no exit flows as the sink node $t$ without a capacity constraint.
3) the decision variables are real numbers, not binary values.
4) each decision variable has an upper bound that is equal to its flow capacity.

The optimization problem is:

Maximize

$$F = \sum_{i=1}^{N} \sum_{j=1}^{N} X_{ij} \tag{57}$$

Subject to:

$$\sum_{j=1}^{N} X_{ji} - \sum_{j=1}^{N} X_{ij} = \begin{cases} -F & if\, i = s \\ F & if\, i = t \\ 0 & otherwise \end{cases} \tag{58}$$

where

$$0 \leq X_{ij} \leq U_k \tag{59}$$

The last column of the PFT contains the upper bound for each of the decision variables. The lower bound of zero is the same for all variables so they need not be in the PFT.

### 4.2.1 Example Problem

This example uses the same cities as the shortest path problem, but the arcs represent railroad service. The source node is a seaport in Los Angeles. The sink node is a transshipment terminal in St. Louis. Table 13 is the PFT that corresponds to this example. The capacities are in units of 100,000 containers per day.

Table 13: PFT for the Maximum Flow Problem Exercise

| Available Links | Network Flow Constraints ||||| Maximize Flow | Capacity $U_k$ |
|---|---|---|---|---|---|---|---|
|  | 1 | 2 | 3 | 4 | 5 |  |  |
| $X_{12}$ | 1 |  |  |  |  | 1 | 5 |
| $X_{13}$ |  | 1 |  |  |  | 1 | 3 |
| $X_{14}$ |  |  | 1 |  |  | 1 | 2 |
| $X_{24}$ | -1 |  | 1 |  |  |  | 5 |
| $X_{25}$ | -1 |  |  | 1 |  |  | 3 |
| $X_{34}$ |  | -1 | 1 |  |  |  | 5 |
| $X_{36}$ |  | -1 |  |  | 1 |  | 3 |
| $X_{45}$ |  |  | -1 | 1 |  |  | 1 |
| $X_{46}$ |  |  | -1 |  | 1 |  | 3 |
| $X_{47}$ |  |  | -1 |  |  |  | 4 |
| $X_{57}$ |  |  |  | -1 |  |  | 5 |
| $X_{67}$ |  |  |  |  | -1 |  | 1 |
| = | 0 | 0 | 0 | 0 | 0 | F |  |

The objective function is the *dot product* of the decision variables and the column next to the last column such that the optimization problem becomes:

Maximize

$$F = X_{12} + X_{13} + X_{14} \tag{60}$$

Subject to:

$$X_{12} - X_{24} - X_{25} = 0 \tag{61}$$

$$X_{13} - X_{34} - X_{36} = 0 \tag{62}$$

$$X_{14} + X_{24} + X_{34} - X_{45} - X_{46} - X_{47} = 0 \tag{63}$$

$$X_{36} + X_{46} - X_{57} = 0 \tag{64}$$

$$X_{36} + X_{46} - X_{67} = 0 \tag{65}$$

where

$$0 \le X_{ij} \le U_k \tag{66}$$

### 4.2.2 Solution Exercise

The solution, illustrated in Figure 17, is that trains can haul 900,000 containers per day from the port in Los Angeles. The arc values in black and red fonts represent the flow volume and the link capacity, respectively. There is currently no constraint on the number of containers that the terminal facilities in St. Louis can receive. However, the problem can include one by adding the constraint that

$$X_{57} + X_{47} + X_{67} \leq U \tag{67}$$

with $U$ being the upper bound capacity.

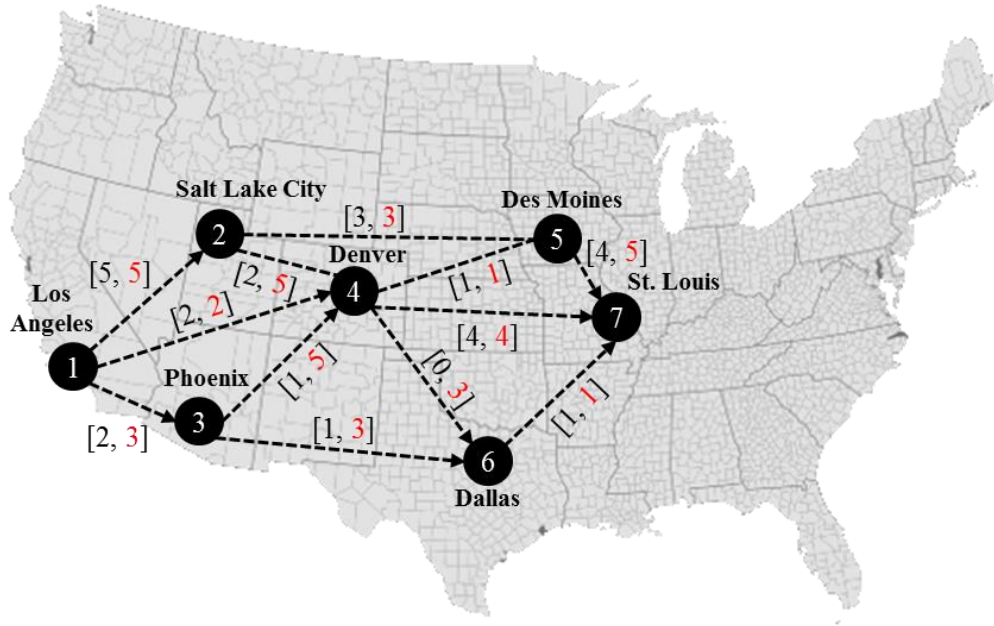

Figure 17: Solution to the example maximum flow problem.

It turns out that this solution is not unique, but the maximum flow is unchanged. Figure 18 illustrates an alternate solution.

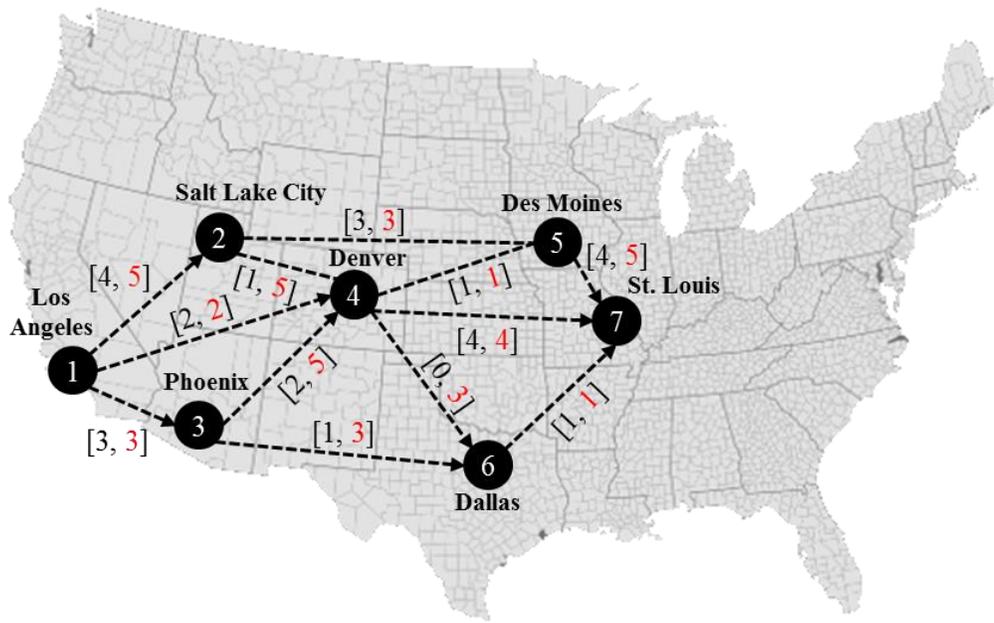

Figure 18: An alternate solution to the example maximum flow problem.

### *4.2.3 Program Display*

```
# Author: Dr. Raj Bridgelall (raj.bridgelall@ndsu.edu)
# Maximum Flow Optimization using scipy.optimize
from IPython import get_ipython
get_ipython().magic('clear')                           # Clear the console
get_ipython().run_line_magic('matplotlib', 'inline')   # plot in the iPython console
import pandas as pd
from pathlib import Path
from scipy.optimize import linprog
import numpy as np
#%%
datapath_in = 'C:/Users/Admin/Documents/Spatial Logistics and Flows/Lab/'
infile = 'Maximum Flow PFT.xlsx'                       # Input filename
filepath_in = Path(datapath_in + infile)               # Path name for untruncated signal
df = pd.DataFrame(pd.read_excel(filepath_in, skiprows=1))   # Problem Formulation
#%%
Nc = df.shape[1]-3      # Number of flow constraints
Nd = df.shape[0]-1      # Number of decision variables
#%%
df = df.fillna(0)                                      # replace all NaN (blanks in Excel) with zeros
A_df = df.iloc[0:Nd,1:Nc+1].astype(int)     # Extract the A parameters
A_parms_eq = A_df.values.transpose()        # Convert the A parametersto a matrix but transposed
b_parms_eq = df.values[-1][1:Nc+1].astype(int)         # Get constraint parameters in an array
#%%
c_parms = df.values.transpose()[-2][:-1].astype(int)   # Get objective parameters in an array
c_parms = c_parms * -1                                 # Transform to a maximization problem
#%%
x_ub = df.values.transpose()[-1][:-1].astype(float)    # List of x upper bounds
x_lb = [0] * len(x_ub)                                 # List of lower bounds (zeros)
#%%
x_bounds = list(zip(x_lb, x_ub))                       # list of lower and upper bound tuples
#%% Call the optimizer
# https://docs.scipy.org/doc/scipy/reference/generated/scipy.optimize.linprog.html
res = linprog(c_parms, A_eq=A_parms_eq, b_eq=b_parms_eq,
              bounds=x_bounds,
              method='revised simplex')
#%% Print the results
print(res); print("\n")
Sol = res.x.astype(int)
print("Values for Decision Variables:\n", Sol)
Arcs = [df.values.transpose()[0][i] for i in range(Nd) if Sol[i] > 0]
Arc_Vals = [Sol[i] for i in range(Nd) if Sol[i] > 0]
print('\n(Arc, Flow):\n', list(zip(Arcs, Arc_Vals)))
print("\nMaximum Flow = ", abs(np.dot(Sol,c_parms)))
```

## 4.3 Warehouse Location Optimization

Some problems require a mix of continuous and binary decision variables. One such problem is the *facility locating problem* where the cost optimization is based on both a per-unit cost and a fixed cost. The per-unit cost includes the costs of producing an item plus transporting the item from facility *i* to consumer *j*. The fixed cost arises if a facility is open and incurs operating, maintenance, rental, and insurance costs.

The general problem is that given *N* facilities and *M* customers, determine which facilities to open and how they could satisfy the demand from all the customers while minimizing cost. The optimization problem is:

Minimize

$$Z = \sum_{i=1}^{N}\sum_{j=1}^{M} c_{ij}Y_{ij} + \sum_{i=1}^{N} f_i X_i \tag{68}$$

Subject to:

$$\sum_{i=1}^{N} Y_{ij} = d_j \quad (j = 1, 2, \ldots, M) \tag{69}$$

and

$$\sum_{j=1}^{M} Y_{ij} \leq u_i X_i \quad (i = 1, 2, \ldots, N) \tag{70}$$

where

$$Y_{ij} \geq 0 \quad (i = 1, 2, \ldots, N),\ (j = 1, 2, \ldots, M) \tag{71}$$

$$X_i \in \{0,1\} \quad (i = 1, 2, \ldots, N) \tag{72}$$

The first part of the objective function is the total variable cost based on the per-unit costs $c_{ij}$ for producing and transporting items from warehouse $i$ to store $j$. The second part of the objective function is the total of the fixed costs $f_i$ for keeping warehouse $i$ open. The decision variable $Y_{ij}$ is the number of units that warehouse $i$ can supply to store $j$ where the demand is $d_j$ units. Hence, the variables are integers with a lower bound of 0. The second constraint by equation (70) assures that the total units shipped from a warehouse is no more than the supply capacity $u_i$ of that warehouse. Hence, if the choice is to close a warehouse, the constraint assures that all $Y_{ij}$ values for that warehouse will be zero.

Table 14 shows the PFT for the warehouse location problem. The last row of the demand constraints contains the number of units of demand $d_j$ for store $j$ for some period such as one day or one week. The last four rows of the cost column contain the fixed cost $f_i$ to keep warehouse $i$ open for the same period. The remaining cells of the cost column stores the per-unit cost for producing and transporting items from warehouse $i$ to store $j$.

Table 14: PFT for the Warehouse Location Problem

| Ship | Stores (Demand Constraint eq) | | | | | Warehouse (Supply Constraint ub) | | | | Cost |
|---|---|---|---|---|---|---|---|---|---|---|
| $i$ to $j$ | 1 | 2 | 3 | 4 | 5 | 1 | 2 | 3 | 4 | ($) |
| $Y_{11}$ | 1 | | | | | 1 | | | | 1 |
| $Y_{12}$ | | 1 | | | | 1 | | | | 2 |
| $Y_{13}$ | | | 1 | | | 1 | | | | 3 |
| $Y_{14}$ | | | | 1 | | 1 | | | | 4 |
| $Y_{15}$ | | | | | 1 | 1 | | | | 5 |
| $Y_{21}$ | 1 | | | | | | 1 | | | 5 |
| $Y_{22}$ | | 1 | | | | | 1 | | | 4 |
| $Y_{23}$ | | | 1 | | | | 1 | | | 3 |
| $Y_{24}$ | | | | 1 | | | 1 | | | 2 |
| $Y_{25}$ | | | | | 1 | | 1 | | | 1 |
| $Y_{31}$ | 1 | | | | | | | 1 | | 1 |
| $Y_{32}$ | | 1 | | | | | | 1 | | 2 |
| $Y_{33}$ | | | 1 | | | | | 1 | | 3 |
| $Y_{34}$ | | | | 1 | | | | 1 | | 4 |
| $Y_{35}$ | | | | | 1 | | | 1 | | 5 |
| $Y_{41}$ | 1 | | | | | | | | 1 | 5 |
| $Y_{42}$ | | 1 | | | | | | | 1 | 4 |
| $Y_{43}$ | | | 1 | | | | | | 1 | 3 |
| $Y_{44}$ | | | | 1 | | | | | 1 | 2 |
| $Y_{45}$ | | | | | 1 | | | | 1 | 1 |
| $X_1$ | | | | | | -60 | | | | 20 |
| $X_2$ | | | | | | | -10 | | | 30 |
| $X_3$ | | | | | | | | -50 | | 20 |
| $X_4$ | | | | | | | | | -55 | 30 |
| $=, \leq$ | 10 | 20 | 30 | 40 | 50 | 0 | 0 | 0 | 0 | Z |

### 4.3.1 Solution Exercise

The solution to satisfy the demand for the example problem is:
Y11 = 5.0
Y12 = 20.0
Y14 = 35.0
Y31 = 5.0
Y33 = 30.0
Y44 = 5.0
Y45 = 50.0
X1 = 1.0
X2 = 0.0
X3 = 1.0
X4 = 1.0

Hence, the solution is to close the second warehouse and achieve a minimum cost of $410 for all suppliers, within the referenced period.

### *4.3.2 Program Display*

```python
from pathlib import Path
from mip import *          # install library from Anaconda prompt: pip install mip
import numpy as np
#%%
datapath_in = 'C:/Users/Admin/DocumentsTL 885/Lectures/10 - Spatial Logistics and Flows/Lab/'
infile = 'Warehouse Location PFT.xlsx'                  # Input filename
filepath_in = Path(datapath_in + infile)                # Path name for untruncated signal
df = pd.DataFrame(pd.read_excel(filepath_in, skiprows=1)) # Read Problem Formulation Table (PFT)
#%%
Ncd = df.iloc[-1,0]         # Number of demand constraints
Ncs = df.shape[1]-2-Ncd     # Number of supply constraints
Nd_C = Ncd * Ncs            # Number of continuous decision variables
Nd_B = Ncs                  # Number of binary decision variables
#%%
df = df.fillna(0)                                       # replace all NaN (blanks in Excel) with zeros
Var_C = df.iloc[0:Nd_C,0]                               # String of continuous variable names
Var_B = df.iloc[Nd_C:-1,0]                              # String of binary variable names
#%%
A_df = df.iloc[0:Nd_C+Nd_B,1:Ncd+Ncs+1]                 # Extract the constraints cols
A_parms_eq = A_df.iloc[:,0:Ncd].astype(int)             # Extract A_eq
A_parms_ub = A_df.iloc[:,Ncd:].astype(float)            # Extract A_ub (already in standard ub form)
b_parms = list(df.values[-1][1:Ncd+1].astype(int))      # Get constraint parameters in an array
c_parms = df.values.transpose()[-1][:-1].astype(int)    # Get objective parameters in an array
#%%
m = Model(solver_name=CBC)                              # use GRB for Gurobi
y = [m.add_var(name=Var_C[i],var_type=INTEGER) for i in range(Nd_C)]   # decision vars y
x = [m.add_var(name=Var_B[i],var_type=BINARY) for i in range(Nd_C,Nd_C+Nd_B)]   # decision vars x
yx = y + x                                              # combined list of decision variables
#%%
m.objective = minimize(xsum(c_parms[i]*yx[i] for i in range(Nd_C+Nd_B))) # add objective function
#%% Add each constraint column
for j in range(Ncd):
    m.add_constr( xsum(A_parms_eq.iloc[i,j]*yx[i] for i in range(len(yx))) == b_parms[j], "ConsY"+str(j) )
for j in range(Ncs):
    m.add_constr( xsum(A_parms_ub.iloc[i,j]*yx[i] for i in range(len(yx))) <= 0, "ConsX"+str(j) )
#%%
Status = m.optimize()
#%% Print the results
print('Model has {} vars, {} constraints and {} nzs'.format(m.num_cols, m.num_rows, m.num_nz))
print('Objective Function: \n',m.objective)
print("Number of Solutions = ", m.num_solutions)
print("Status = ", Status)
#%%
ArcsY = [str(m.vars[i]) for i in range(Nd_C) if y[i].x != 0]
ArcYVals = [y[i].x for i in range(Nd_C) if y[i].x != 0]
for i in range(len(ArcsY)):
    print('{} = {} '.format(ArcsY[i], round(ArcYVals[i],3)))
#%%
VarsX = [str(m.vars[i]) for i in range(Nd_C,Nd_C+Nd_B) ]
Xvals = [x[i].x for i in range(Nd_B) ]
for i in range(Nd_B):
    print('{} = {} '.format(VarsX[i], Xvals[i]))
#%%
print('Minimum Cost = ${}'.format(m.objective_value))
print('Constraints:')
for j in range(m.num_rows):
    print(m.constrs[j])
```

# Chapter 5 - Summary and Conclusions

This tutorial introduced laboratory exercises in linear programming that featured several important optimization problems in supply chain management and transport logistics. Chapter 1 introduced the various types of integer programming and presented a cognitive framework to formulate the optimization problem by organizing the facts of the problem in a standard manner. The problem formulation table (PFT) helped to identify all the decision variables, their linear relationship with the constraints and objective function, and their value bounds. The PTF becomes impractical for importing an exceptionally large and mostly sparse table for problems that have many constraints and decision variables. Therefore, some of the exercises show how to use the multiple integer programming (MIP) library in Python to scale the PFT for larger problems. Chapter 1 described how to integrate two free GIS tools into the optimization workflow to visualize solutions to the optimization problems. The remaining chapters cover optimization problems involving mobility, spatial placements, and distributions in a supply chain. The mobility optimization problems determined the shortest path in a network and the minimum cost tour to visit all nodes once at most. The spatial optimization problems explained neighborhood coverage with the least number of facilities, capturing the maximum flows in a network, creating non-competing zones, and optimizing service locations. The spatial logistics problems demonstrated distribution optimizations to meet demands within capacity constraints at the minimum cost, maximizing flows through a capacity-constrained network, and locating warehouses to minimize variable costs and fixed operating costs. Students and practitioners can modify the demonstration code and GIS workflow to solve their own optimization problems.